\documentclass[12pt,reqno]{amsart}
\usepackage{amsmath , amssymb }
\usepackage{graphics}

\usepackage[mathscr]{eucal}

\newdimen\unit\newdimen\psep\newcount\nd\newcount\ndx\newbox\dotb\newbox\ptbox
\newdimen\dx\newdimen\dy\newdimen\dxx\newdimen\dyy\newdimen\hgt
\newdimen\xoff\newdimen\yoff
\newcommand\clap[1]{\hbox to 0pt{\hss{#1}\hss}}
\newcommand\vdisk[1]{{\font\dotf=cmr10 scaled #1\dotf.}}
\newcommand\varline[2]{\setbox\dotb\hbox{\vdisk{#1}}\xoff=-.5\wd\dotb
\wd\dotb=0pt\yoff=-.5\ht\dotb\psep=#2\ht\dotb}
\newcommand\varpt[1]{\setbox\ptbox\clap{\vdisk{#1}}\setbox\ptbox
\hbox{\raise-.5\ht\ptbox\box\ptbox}}
\newcommand\cpt{\copy\ptbox}
\newcommand\point[3]{\rlap{\kern#1\unit\raise#2\unit\hbox{#3}}}
\newcommand\setnd[4]{\dx=#3\unit\advance\dx-#1\unit\divide\dx by\psep
\dy=#4\unit\advance\dy-#2\unit\divide\dy by\psep \multiply\dx
by\dx\multiply\dy by\dy\advance\dx\dy\nd=1\advance\dx-1sp
\loop\ifnum\dx>0\advance\dx-\nd sp\advance\nd1\advance\dx-\nd
sp\repeat}

\newcommand\dline[5]{{\nd=#5\hgt=#2\unit\dx=#3\unit\advance\dx-#1\unit
\divide\dx by\nd\dy=#4\unit\advance\dy-#2\unit\divide\dy by\nd
\advance\hgt\yoff\rlap{\kern#1\unit\kern\xoff\loop\ifnum\nd>1\advance\nd-1
\advance\hgt\dy\kern\dx\raise\hgt\copy\dotb\repeat}}}

\newcommand\qellip[4]{{\setnd{0}{0}{#3}{#4}\dx=\unit\dy=0pt\raise\yoff\rlap{%
\kern#1\unit\kern\xoff\raise#2\unit\hbox{\loop\ifnum\dx>0\rlap{\kern#3\dx
\raise#4\dy\copy\dotb}\hgt=\dx\divide\hgt
by\nd\advance\dy\hgt\hgt=\dy \divide\hgt
by\nd\advance\dx-\hgt\repeat\rlap{\raise#4\dy\copy\dotb}}}}}
\newcommand\bez[6]{{\setnd{#1}{#2}{#3}{#4}\ndx=\nd\setnd{#3}{#4}{#5}{#6}
\ifnum\ndx>\nd\nd=\ndx\fi\dx=#3\unit\advance\dx-#1\unit\dy=#4\unit
\advance\dy-#2\unit\dxx=#5\unit\advance\dxx-#1\unit\dyy=#6\unit\advance
\dyy-#2\unit\advance\dxx-2\dx\advance\dyy-2\dy\divide\dxx
by\nd\divide\dyy by\nd\advance\dx.25\dxx\advance\dy.25\dyy\divide\dx
by\nd\divide\dy by\nd \multiply\nd
by2\dx=100\dx\dy=100\dy\dxx=100\dxx\dyy=100\dyy\divide\dxx by\nd
\divide\dyy by\nd\hgt=#2\unit\raise\yoff\rlap{\kern#1\unit\kern\xoff
\raise\hgt\copy\dotb\loop\ifnum\nd>0\advance\nd-1\advance\hgt0.01\dy
\kern0.01\dx\raise\hgt\copy\dotb\advance\dx\dxx\advance\dy\dyy\repeat}}}
\newcommand\ptu[3]{\point{#1}{#2}{\cpt\raise1ex\clap{$\scriptstyle{#3}$}}}
\newcommand\ptd[3]{\point{#1}{#2}{\cpt\raise-1.8ex\clap{$\scriptstyle{#3}$}}}
\newcommand\ptr[3]{\point{#1}{#2}{\cpt\raise-.4ex\rlap{$\ \scriptstyle{#3}$}}}
\newcommand\ptl[3]{\point{#1}{#2}{\cpt\raise-.4ex\llap{$\scriptstyle{#3}\ $}}}
\newcommand\ptlu[3]{\point{#1}{#2}{\raise.8ex\clap{$\scriptstyle{#3}$}}}
\newcommand\ptld[3]{\point{#1}{#2}{\raise-1.6ex\clap{$\scriptstyle{#3}$}}}
\newcommand\ptlr[3]{\point{#1}{#2}{\raise-.4ex\rlap{$\,\scriptstyle{#3}$}}}
\newcommand\ptll[3]{\point{#1}{#2}{\raise-.4ex\llap{$\scriptstyle{#3}\,$}}}

\newcommand\thnline{\varline{400}{.6}}

\varpt{2500}\thnline\unit=2em

\newtheorem{thm}{Theorem}
\newtheorem*{Bary}{Lemma (Baranyai, 1973)}

\newtheorem{conj}{Conjecture}
\newtheorem{prop}[thm]{Proposition}
\newtheorem{qu}{Question}
\newtheorem{prob}{Problem}
\newtheorem{lemma}[thm]{Lemma}
\newtheorem{cor}[thm]{Corollary}

\newtheorem{obs}[thm]{Observation}
\theoremstyle{definition}\newtheorem{rmk}{Remark}
\theoremstyle{definition}
\theoremstyle{definition}\newtheorem{eg}{Example}
\newcommand{\ds}{\displaystyle}
\newcommand{\ul}{\underline}

\def\Ex{\mathbb{E}}

\def\N{\mathbb{N}}

\def\RR{\mathbb{R}}
\def\Z{\mathbb{Z}}

\def\Pr{\mathbb{P}}
\def\le{\leqslant}
\def\ge{\geqslant}

\def\eps{\varepsilon}

\begin{document}
\title[Majority bootstrap on the hypercube]{Majority bootstrap percolation on the hypercube}

\author{J\'ozsef Balogh}
\address{Department of Mathematics\\ University of Illinois\\ 1409 W. Green Street\\ Urbana, IL 61801, USA} \email{jobal@math.uiuc.edu}

\author{B\'ela Bollob\'as}
\address{Trinity College\\ Cambridge CB2 1TQ\\ England\\ and \\ Department of Mathematical Sciences\\ The University of Memphis\\ Memphis, TN 38152, USA} \email{B.Bollobas@dpmms.cam.ac.uk}

\author{Robert Morris}
\address{Instituto Nacional de Matem\'atica Pura e Aplicada, Estrada Dona Castorina, 110,
 Jardim Bot\^{a}nico, Rio de Janeiro, Brazil} \email{rdmorrs1@impa.br}\thanks{The first author was supported during this research by OTKA grant T049398 and NSF grants DMS-0302804,  DMS-0603769 and DMS 0600303, and UIUC Campus Research Board 06139 and 07048, the second by ITR grant CCR-0225610 and ARO grant W911NF-06-1-0076, and the third by PRONEX CNPq/FAPERJ grant E-26/171.167/2003-APQ1 and by MCT grant PCI EV-8C}

\begin{abstract}
In majority bootstrap percolation on a graph $G$, an infection spreads according to the following deterministic rule: if at least half of the neighbours of a vertex $v$ are already infected, then $v$ is also infected, and infected vertices remain infected forever. Percolation occurs if eventually every vertex is infected.

The elements of the set of initially infected vertices, $A \subset V(G)$, are normally chosen independently at random, each with probability $p$, say. This process has been extensively studied on the sequence of torus graphs $[n]^d$, for $n = 1,2, \ldots$, where $d = d(n)$ is either fixed or a very slowly growing function of $n$. For example, Cerf and Manzo~\cite{CM} showed that the critical probability is $o(1)$ if $d(n) \le \log_{*} n$, i.e., if $p = p(n)$ is bounded away from zero then the probability of percolation on $[n]^d$ tends to one as $n \to \infty$. 

In this paper we study the case when the growth of $d$ to $\infty$ is not excessively slow; in particular, we show that the critical probability is $1/2 + o(1)$ if $d \ge (\log \log n)^2 \log \log \log n$, and give much stronger bounds in the case that $G$ is the hypercube, $[2]^d$.
\end{abstract}

\maketitle

\section{Introduction}\label{intro}

Consider a finite graph $G$ with two parameters, $q_v$ and $r_v$, attached to each vertex $v$, with $q_v + r_v$ greater than the degree of $v$. Suppose each of the vertices may take either one of two states, `active' and `inactive', say. At each instant, some of the vertices may `wake up' at random; whenever a vertex $v$ does so, if at least $q_v$ of its neighbours are active then it becomes active, if at least $r_v$ of its neighbours are inactive then it becomes inactive, and if neither of these cases holds then $v$ keeps its state. Given an initial distribution of active sites, one can then ask what happens to the system in the long run.

For example, we may take $G$ to be a $d$-regular graph with $d$ odd, and parameters $q_v = r_v = (d+1)/2$ for every vertex $v$: in this case, whenever a vertex wakes up, its state becomes that of the majority of its neighbours. For a $d$-regular graph with $d$ even, we may take $q_v = r_v = d/2 + 1$ for every vertex $v$.

A special (and very well-studied) example of this process is the zero-temperature Ising model, where the process occurs on the lattice $\Z^d$. Despite all this interest, however, only rather weak (though far from easy) results have so far been proven (see \cite{FSS}, \cite{NNS} and \cite{NS}, for example). Another example is the celebrated model of the brain, introduced over 60 years ago by McCulloch and Pitts~\cite{MP}, and shown by them to be complex enough to include a universal Turing machine as a particular case.

In this paper we shall study a simpler model, introduced by Chalupa, Leith and Reich~\cite{CLR} in 1979, in which the vertices may only change state in one direction, from inactive to active, say. In order to make this easier to remember, let us refer instead to `healthy' and `infected' vertices, so that a healthy vertex may be infected, but infected vertices never recover.

More generally (and more precisely), let $G$ be a finite graph, let $r \in \N$, and let $A \subset V(G)$ be a set of initially infected vertices. The set of infected vertices is updated as follows: if a healthy vertex has at least $r$ infected neighbours, he becomes infected; otherwise he remains healthy. In other words, we have a sequence of sets
$$A = A^{(0)} \subset A^{(1)} \subset \ldots \subset A^{(m)} \subset \ldots$$
where $A^{(m+1)} = A^{(m)} \cup \{v \in V(G) : |\Gamma(v) \cap A^{(m)}| \ge r\}$. If the entire graph is eventually infected, i.e., $A^{(m)} = V(G)$ for some $m \in \N$, then the set $A$ is said to \emph{percolate} on $G$. This process is known as the \emph{$r$-neighbour bootstrap percolation on $G$}; if $G$ is $d$-regular with $d$ odd and $r = \lceil d/2 \rceil$, then this is the \emph{majority} bootstrap percolation. We remark that this process is very different from the \emph{random majority} bootstrap process studied in~\cite{BBJW}, and also the \emph{biased majority} process studied in~\cite{Schbias}, since in each of those processes vertices may change states in both directions, i.e., may be infected and then later healed.

The bootstrap process has been well-studied in the case that $G$ is the torus $[n]^d$, where $r$ and $d$ are both fixed (with $r \le d$), the elements of the set $A$ are chosen independently at random, and $n \to \infty$. (In fact the results below were proved for the grid (i.e., the subgraph of $\Z^d$ induced by the vertices of $[n]^d$) but the proofs also apply to the corresponding tori. Note that for $n = 2$ the grid is $d$-regular, but for $n \ge 3$ the degrees vary from $d$ to $2d$.)

Let $p = p(n) = \Pr(v \in A)$ for each $v \in V(G)$, and write $\Pr_p$ for the corresponding (product) probability measure. Clearly the probability that $A$ percolates is monotone in $p$, since extra initial infections can only make percolation more likely. Hence there exists a (unique) value $p_c \in (0,1)$, depending on $n$, $d$ and $r$, for which $\Pr(A$ percolates$) = 1/2$. We call this value the \emph{critical probability}, and in general define
$$p_c(G,r) \; = \; \sup\left\{p \in (0,1)\, : \, \Pr_p\left(A\textup{ percolates on }G \right) \le \frac{1}{2}  \right\}.$$
The challenge is to determine the value of $p_c$, and also the size of the \emph{critical window}, i.e., the range $p_{1-\eps} - p_\eps$, where $p_\alpha$ satisfies $\Pr_{p_\alpha}(A$ percolates$) = \alpha$ for each $\alpha \in (0,1)$, and $\eps \to 0$. In general we would like to show that the window is small, and hence that the threshold for percolation is `sharp'.

The first rigorous results for bootstrap percolation on finite graphs were obtained by Aizenman and Lebowitz~\cite{AL}, who showed that
$$\frac{c(d)}{(\log n)^{d-1}} \; \le \; p_c([n]^d,2) \; \le \; \frac{C(d)}{(\log n)^{d-1}}$$ for some functions $0 < c(d) < C(d)$, and moreover that the size of the critical window for $[n]^d$ is $O(p_c)$ when $d$ is fixed and $r = 2$. The problem when $3 \le r \le d$ seems to be more difficult (when $r > d$ we have $p_c = 1 - o(1)$, so the problem is less interesting). Despite this, Cerf and Manzo~\cite{CM}, building on work of Cerf and Cirillo~\cite{CC} (as well as the older work of Schonmann~\cite{Sch92} on the lattice $\Z^d$), were able to prove the corresponding result for all fixed $d$ and $r$. They proved that, if $r \le d$, then $$\frac{c(d,r)}{(\log_{r-1} n)^{d-r+1}} \; \le \; p_c([n]^d,r) \; \le \; \frac{C(d,r)}{(\log_{r-1} n)^{d-r+1}},$$ for some functions $0 < c(d,r) < C(d,r)$, where $\log_{r-1} n$ is the $(r-1)$-times iterated logarithm, i.e., $\log_1 n = \log n$ and $\log_{k+1} n = \log (\log_k n)$ for each $k \in \N$. In particular, their proof implies that $p_c([n]^d,d) = o(1)$ if $1 \ll d \le \log_{*} n$, where $\log_{*} n$ is the number $k$ such that $\log_k n \ge 1 > \log_{k+1} n$. In particular, note that $\log_{*} n \ll \log_k n$ for every $k \in \N$.

Finally, we remark that in the case $d = r = 2$, Holroyd~\cite{Hol} was able to prove a much sharper result: that in fact
$$p_c([n]^2,2) = \ds\frac{\pi^2}{18\log n} + o\left(\ds\frac{1}{\log n}\right).$$
(Here, and throughout, $\log$ is taken to the base $e$). The reader who is interested in bootstrap percolation on other types of graphs should also see the work of Balogh and Pittel~\cite{BP} on random $d$-regular graphs, and of Balogh, Peres and Pete~\cite{BPP} on infinite trees. An application of the techniques of bootstrap percolation to the zero-temperature Ising model may be found in~\cite{FSS}, and for a brief survey of the physical applications of the bootstrap process, see~\cite{braz}.

In this paper we shall be interested in majority bootstrap percolation on very high-dimensional lattices, and more generally on arbitrary $d$-regular `lattice-like' graphs. The fundamental result of Cerf and Manzo, stated above, has two drawbacks.  The first is that the known bounds on $c(d,r)$ and $C(d,r)$ are rather far apart for large $d$ and $r$, and the second that the theorem is only useful when $n$ is extremely large (so that $\log_{r-1} n > c(d,r)$). We shall show that the latter problem is unavoidable; indeed when $n$ is not quite so large (at most $2^{2^{\sqrt{\frac{d}{\log d}}}}$, say) then the behaviour of $p_c$ is quite different (for sufficiently large $d$). We shall also study in more detail majority bootstrap percolation on one particular graph, the hypercube $Q_d = [2]^d$, which may be thought of as an extreme point of the family $\{[n]^d\}$, where $n = n(t)$ and $d = d(t)$ are arbitrary functions. The hypercube is a very well-studied combinatorial object; for example, see the work (relating to a different percolation problem on the hypercube) of~\cite{ES}, \cite{AKS} and \cite{BKL}, and the recent improvements of~\cite{BCHSS}, \cite{HS1} and \cite{HS2}. In a subsequent paper~\cite{BBMn^d} we shall investigate $2$-neighbour bootstrap percolation on high-dimensional lattices, and show that the problem essentially reduces to the equivalent question for the hypercube. As the reader will discover, a similar phenomenon also occurs for majority percolation.

The structure of the remainder of the paper is as follows. In the next section we shall state our main results; in Section~\ref{tools} we shall describe some of the simple tools that we shall use later; in Sections~\ref{uppersec} and \ref{lowersec} we shall prove fairly strong bounds on $p_c(Q_n,n/2)$; in Section~\ref{dregsec} we shall study more general $d$-regular graphs; and in Section~\ref{qusec} we shall describe some open problems and conjectures. The paper ends with an appendix in which the tools from Section~\ref{tools} are proved.

\section{Main Results}

In this section we shall state our main results. We begin with a theorem which gives bounds on the critical probability for majority bootstrap percolation on the hypercube. Recall that the $n$-dimensional hypercube, $Q_n = [2]^n$, is the $n$-regular graph with vertex set $V(Q_n) = \{0,1\}^n$, and edge set $E(Q_n)$, where $$xy \in E(Q_n)\textup{ if and only if }|\{i : x_i \neq y_i\}| = 1.$$ Observe that the vertices of $Q_n$ may also be thought of as subsets of $[n]$, in which case $xy \in E(Q_n)$ if and only if $|x \triangle y| = 1$.

We shall assume throughout that $n$ is even, so that we may write $n/2$ instead of $\lceil n/2 \rceil$. However, all our proofs are also valid for $n$ odd.

\begin{thm}\label{sharp}
Let $n \in \N$, $\lambda \in \RR$, $$p \; = \; p(n) \; = \; \frac{1}{2} \: - \: \frac{1}{2} \sqrt{\ds\frac{\log n}{n}} \: + \: \frac{\lambda \log \log n}{\sqrt{n \log n}},$$
and let the elements of $A \subset V(Q_n)$ be chosen independently at random, each with probability $p$.
Then in majority bootstrap percolation, $$\Pr(A\textup{ percolates on }Q_n) \to \left\{
\begin{array} {r@{\quad \textup{if} \quad}l} 0 & \lambda \le -2 \\[+0.5ex]
1 & \lambda > 1/2
\end{array}\right.$$ as $n \to \infty$.
In particular,
\begin{align*}
& \hspace{1cm} p_c\big( Q_n, n/2 \big) \; \ge \;
\ds\frac{1}{2} \: - \: \frac{1}{2} \sqrt{\ds\frac{\log n}{n}} \: - \: \ds\frac{ 2\log \log n }{ \sqrt{n \log n} }\\
& \textup{\emph{for sufficiently large }}n\textup{\emph{, and}} \\
& \hspace{1cm} p_c\big( Q_n, n/2 \big) \; \le \;
\ds\frac{1}{2} \: - \: \frac{1}{2} \sqrt{\ds\frac{\log n}{n}} \: + \: \ds\frac{ \log \log n }{ 2\sqrt{n \log n} } + o\left(\ds\frac{ \log \log n}{2\sqrt{n \log n} } \right)
\end{align*}
as $n \to \infty$.
\end{thm}

Note that Theorem~\ref{sharp} determines the first two terms in the expansion of $p_c\big(Q_n,n/2\big)$, but not even the order of the third, since it does not tell us whether or not $A$ is likely to percolate when $\lambda = 0$. In Section~\ref{lowersec}, we shall discuss how these bounds might be further improved, and explain the limitations of our method.

We next turn to more general $d$-regular graphs on $N$ vertices. As it turns out, the method of the proof of Theorem~\ref{sharp} can be adapted, in a slightly simpler form, to deal with a wide range of such graphs. As usual, we write $\Gamma(u)$ for the set of neighbours of a vertex $u \in V(G)$. Furthermore, for each $u \in V(G)$ and $k \in \N$ let us define
$$S(u,k) \; = \; \{v \in V(G)\, : \, d(u,v) = k\},$$ so that $S(u,1) = \Gamma(u)$, and also let $B(u,k) = \{v : d(u,v) \le k\}$. The property of the hypercube which we shall need in order to prove the theorem below, is that the set $S(u,k) \cap \Gamma(v)$ is not `too big' for any $u \in V(G)$, $v \in V(G) \setminus B(u,k-1)$, and `sufficiently large' $k$ (the exact size needed depends on $N$).

The following theorem is somewhat technical, so on a first reading the reader may wish to assume that $k$ and the functions $f_i$ are all constant.

\begin{thm}\label{dreg}
For each $d \in \N$, let $N = N(d),k = k(d) \in \N$, and let $G = G(d)$ be a $d$-regular graph on $N$ vertices. Furthermore, let $\omega, f, f_1, \ldots, f_k: \N \to \N$ be arbitrary functions satisfying
$$1 \, \le \, f_i(d) \, \le \, f(d) \, = \, o\left(\ds\frac{d}{k\log d} \right)$$ and $\omega(d) \to \infty$ as $d \to \infty$. Suppose that, for every $d \in \N$ and $i \in [k]$
$$|S(x,i) \cap \Gamma(y)| \; \le \; f_i(d)$$
for every $x \in V(G)$ and $y \in V(G) \setminus B(x,i-1)$, and
$$N \; \le \; \exp\left( \frac{d^k}{\big(\omega(d)k\big)^k \big(f_{k-1}(d) + f_k(d)\big) \prod_{i=1}^{k-1} f_i(d)} \right)$$ for every $d \in \N$.
Then
$$p_c\big(G,d/2\big) \, = \, \frac{1}{2} + o(1)$$ as $d \to \infty$.
\end{thm}

We remark that the bound on $N$ in Theorem~\ref{dreg} cannot be improved substantially. An example showing this will be described in Section~\ref{dregsec}. Other extensions are possible though; in particular, the bound on $|S(x,i) \cap \Gamma(y)|$ does not have to hold for every $x \in V(G)$: rather, a (small) exceptional set is permissible, since with high probability these vertices will be infected in the first round anyway. We shall not need this generalization however, and do not wish to further complicate the statement of the theorem unnecessarily.

Our main motivation for proving Theorem~\ref{dreg} is the following, almost immediate corollary. Let $[n]^d$ denote the $d$-dimensional $n \times \ldots \times n$ torus, and note that $[n]^d$ is $2d$-regular, and has $n^d$ vertices. Moreover, $[n]^d$ satisfies the conditions of Theorem~\ref{dreg} with $f_i(d) = i + 1$ for each $i \in \N$ (see the `Proof of Corollary~\ref{n^d}' in Section~\ref{dregsec}). Thus, we may apply the theorem as long as $k^2 = o\left(\ds\frac{d}{\log d} \right)$, and $n^d \le 2^{2^k}$; doing so gives the following result.

\begin{cor}\label{n^d}
Let $n = n(t)$ and $d = d(t)$ be functions satisfying
$$3 \; \le \; n \; = \; 2^{2^{O\left(\sqrt{\frac{d}{\log d}}\right)}},$$
or equivalently, $d \ge \eps(\log \log n)^2\log\log\log n$ for some $\eps > 0$. Then
$$p_c\big( [n]^d, d \big) \; = \; \frac{1}{2} \, + \, o(1)$$ as $t \to \infty$.
\end{cor}

The proofs of Theorems~\ref{sharp} and \ref{dreg} are not short, but the main ideas may be summarised in a few sentences. For the upper bound in Theorem~\ref{sharp}, we study the first two steps of the process in detail, and show that with high probability, at least (about) $3/4$ of the vertices will be infected by this stage; the rest of the proof is then simply a matter of battling the weak dependence between the events $\{x \in A^{(2)}\}_{x \in V(Q_n)}$.

The upper bound is somewhat harder, since the process may continue running for many steps. We overcome this by introducing a new, `more generous' process, which nonetheless stops quickly; a simple coupling then shows that the original process must also stop eventually, before it has infected the entire vertex set. The proof of Theorem~\ref{dreg} is similar, though in this case the process must be allowed to run for many more steps, since we are dealing with much larger vertex sets. However, since we only wish to prove a much weaker result (that $p_c(G,d) = 1/2 + o(1)$), the details of the calculations become much simpler.

\section{Tools}\label{tools}

Since all the results in this section will either be well known, or simple approximations of binomials, we shall postpone the proofs until the appendix. We begin by recalling the standard Chernoff bound (see \cite{RG} for example).

\begin{lemma}\label{normalcher}
Let $n \in \N$, $0 < p < 1$, $t \ge 0$ and $S(n) \sim \textup{Bin}(n,p)$. Then
$$(a) \hspace{2cm} \Pr\big(S(n) \ge np + t\big) \; \le \; \exp\left( -\ds\frac{2t^2}{n} \right)\hspace{2cm}$$
and similarly,
$$(b) \hspace{2cm} \Pr\big(S(n) \le np - t\big) \; \le \; \exp\left( -\ds\frac{2t^2}{n} \right)\hspace{2cm}$$
\end{lemma}

The following lemma gives an almost matching lower bound in the case in which we shall be interested.

\begin{lemma}\label{chernoff}
Let $C \ge 0$ be a constant, and $n \in \N$ be sufficiently large. Let $p = \ds\frac{1}{2} - \delta$, where $0 \le 8\delta^4 n \le 1$, and let $S(n) \sim B(n,p)$. Then
$$\Pr\left(S(n) \ge \frac{n}{2} + C\right) \; \ge \;\exp\left( - 2\delta^2 n \, - \, 4\delta \sqrt{\frac{n}{\log n}} \, - \, \frac{\log \log n}{2} \, - \, 6\right).$$
\end{lemma}

We next state a simple generalization of Lemma~\ref{normalcher}, which will be a key tool in the proofs of the lower bounds in Theorems~\ref{sharp} and~\ref{dreg}.

\begin{lemma}\label{layer4}
Let $t,k,d_1,\ldots,d_k \in \N$ and $p \in (0,1)$. Let $X_i \sim \textup{Bin}(d_i,p)$ for each $i \in [k]$, let $Y_k = \sum_{i=1}^k iX_i$, and let $D(k) = \sum_{i=1}^k i^2 d_i$. Then
$$\Pr\big(Y_k \ge \Ex(Y_k) + t\big) \; \le \; (2t)^{k-1}\exp\left( - \frac{2t^2}{D(k)} \right).$$
\end{lemma}

Now an easy approximation, will we shall use to prove the upper bounds in Theorems~\ref{sharp} and \ref{dreg}.

\begin{lemma}\label{nunlikely}
Let $p \in (0,1)$ and $n \in \N$ satisfy $pn^2 \le 1$, and let $S(n) \sim \textup{Bin}(n,p)$. Then
$$\Pr\big(S(n) \ge m\big) \; \le \; 2p^{m/2}$$ for every $m \in [n]$. In particular, if $c, \eps > 0$ and $p \le e^{-c n}$, then, for some $b = b(c,\eps) > 0$ not depending on $n$,
$$\Pr(S(n) \ge \eps n) \; \le \; e^{-bn^2}.$$
\end{lemma}

In order to apply Lemma~\ref{nunlikely}, we shall make frequent use of the following simple lemma.

\begin{lemma}\label{partition}
Let $G$ be a graph, let $k,m \in \N$, and suppose that for each $x \in V(G)$, $$|B(x,k)| \: = \: |\{y \in V(G) : d(x,y) \le k\}| \: \le \: m.$$
Then there exists a partition $$V(G) \: = \: B_1 \cup \ldots \cup B_{m}$$ of $V(G)$, such that if $y,z \in B_i$, then $d(y,z) \ge k+1$.
\end{lemma}

Lemma~\ref{partition} immediately implies the following result for hypercubes. We remark that below, and throughout the paper, we shall often write $x \in Q_n$ to mean $x \in V(Q_n)$.

\begin{lemma}\label{hyperpartition}
Let $n,k \in \N$, and $x \in Q_n$. Then there exists a partition $$S(x,k) \: = \: B_1 \cup \ldots \cup B_{m}$$ of $S(x,k)$ into $m \le k{n \choose {k-1}} \le 2n^{k-1}$ sets, such that if $x,y \in B_j$ for some $j$, then $d(x,y) \ge 2k$.
\end{lemma}

We remark that Lemma~\ref{hyperpartition} is a simple special case of the following old result, due to Baranyai~\cite{Bary}.

\begin{Bary}
Let $K^{(h)}_n$ denote the complete $h$-uniform hypergraph on $n$ vertices. If $h$ divides $n$, then $K^{(h)}_n$ can be partitioned into $1$-regular hypergraphs.
\end{Bary}

We shall also need the following two, rather easy lemmas, and one trivial observation, which follows by the convexity of ${x \choose 2}$.

\begin{lemma}\label{giveone}
Let $\delta = \delta(n) \to 0$ as $n \to \infty$, and let $p = \ds\frac{1}{2} - \delta$. Let $S(n) \sim
B(n,p)$ and let $S'(n) \sim 1 + B(n-1,p)$. Then, for any $0 \le m = m(n) \le n/2$,
$$\Pr\left(S'(n) \ge m \right) \; = \; \big( 1 + o(1) \big)\,\Pr\left(S(n) \ge m \right)$$ as $n \to \infty$.
\end{lemma}

\begin{lemma}\label{silly}
Let $\delta = \delta(n) \to 0$ as $n \to \infty$, and let $p = \ds\frac{1}{2} - \delta$. Let $S(n) = X + Y(n)$, where $X \sim \textup{Bin}(1,p)$ and $Y(n) \sim \textup{Bin}(n-1,p)$. Then, for any $0 \le m = m(n) \le n/2$,
$$\Pr\left(X(1) = 1 \,|\, S(n) \ge m\right) \; = \; \big(1 + o(1)\big)\Pr\big(X(1) = 1\big)$$ as $n \to \infty$.
\end{lemma}

\begin{obs}\label{convexity}
Let $k, a_1, \ldots, a_k, A \in \N$, and suppose $\ds\max_i\{a_i\} \le A$. Then
$$\sum_i {a_i \choose 2} \; \le \; \frac{\sum_i a_i}{A} {A \choose 2}.$$
\end{obs}

Finally, we shall use the following, probably well-known lemma (see \cite{LP}, Lemma B.7), and we recall the standard second moment method (see \cite{RG}, for example).

\begin{lemma}\label{morethanhalf}
Let $p \in (0,1)$, let $n \in \N$, and let $S(n) \sim \textup{Bin}(n,p)$. Then
$$\Pr\left(S(n) \le \lfloor np \rfloor - 1\right) \; \le 1/2 \; \le \; \Pr\left(S(n) \le \lceil np \rceil\right).$$
In particular, if $\eps > 0$ and $p = p(n) \in [\eps,1-\eps]$ for every $n \in \N$, then
$$\Pr\left(S(n) \ge np\right) = 1/2 + o(1)$$
as $n \to \infty$.
\end{lemma}

\begin{lemma}\label{2mm}
For any random variable $X$, and any $\alpha > 0$, we have
$$\Pr\Big(|X - \Ex(X)| \: \ge \: \alpha |\Ex(X)|\Big) \; \le \; \frac{\textup{Var}(X)}{\alpha^2 \Ex(X)^2}.$$ In particular, if $X_n$ is a sequence of non-negative random variables, such that $\textup{Var}(X_n) = o(\Ex(X_n)^2)$, then
$$\Pr\left(\frac{\Ex(X_n)}{2} \: \le \: X_n \: \le \: 2\Ex(X_n)\right) \: = \: 1 - o(1).$$
\end{lemma}

\section{Proof of the upper bound in Theorem~\ref{sharp}}\label{uppersec}

Let $n \in \N$ be sufficiently large, and let $$p(n) \; = \; \frac{1}{2} \: - \: \frac{1}{2}\sqrt{\frac{\log n}{n}} \: + \: \frac{\lambda \log \log n}{\sqrt{n \log n}},$$ for some $\lambda \in \RR$. We shall show that if $\lambda > 1/2$, then for some constant $c > 0$ and each vertex $x \in V$, $\Pr(x \not\in A^{(11)}) < e^{-cn^2}$, and thus that $\Pr(A^{(11)} = V) = 1 - o(1)$.

The proof comes in three stages: first we shall show that $\Pr(x \in A^{(2)}) > 3/4$ (Lemma~\ref{3/4}); next we shall show that $\Pr(x \notin A^{(5)}) < e^{-cn}$ for some constant $c$ (Lemma~\ref{2to5}); and finally, we shall show that $\Pr(x \not\in A^{(11)}) < e^{-cn^2}$ (Lemma~\ref{exptoall}).

In all that follows, we shall assume that the elements of $A^{(0)}$ are chosen independently, each with probability $p$, and that $A^{(i)}$ are the infected vertices after $i$ rounds.

In fact, we need a slightly more general concept. Let $A_r^{(i)}$ denote the set of infected vertices after round $i$ if the infection threshold is $r$. Thus $A^{(i)} = A^{(i)}_{n/2}$.

\begin{lemma}\label{3/4}
Let $x \in V(Q_n)$, and $p$ and $A^{(0)}$ be defined as above. Let $r \le n/2 + 100$. If $\lambda > 1/2$, then $\Pr(x \in A_r^{(2)}) \ge 3/4 + o(1)$ as $n \to \infty$.
\end{lemma}

\begin{proof}
Let $x \in V(Q_n)$, where $n$ is large. With probability $p = 1/2 + o(1)$, we have $x \in A^{(0)}$. For those $x \notin A^{(0)}$, we shall use Lemma~\ref{chernoff} to show that $\Pr\big(|\Gamma(x) \cap A_r^{(1)}| \ge r \big) \ge 1/2 + o(1)$.

So let us assume that $x \notin A^{(0)}$, and let $R = \Gamma(x) \cap A^{(0)}$ and $S = \Gamma(x) \cap A_r^{(1)} \setminus A^{(0)}$. We want to show that $\Pr(|R| + |S| \ge r) \ge 1/2 + o(1)$. This follows easily from the following claim.\\[-1ex]

\noindent\ul{Claim}: Let $100 \le m \le n/2$. Then $$\Pr\left( |S| \le \ds\frac{\sqrt{n} (\log n)^{(4\lambda - 1)/2}}{e^{11}} \; \Big| \; |R| \ge r - m \right) = o(1).$$

\begin{proof}[Proof of claim]
We use the standard second moment method (Lemma~\ref{2mm}). Throughout the proof of the claim we assume that $|R| \ge r - m$; in particular, we shall sometimes write $\Pr(\cdot)$ for $\Pr(\cdot \, | \, |R| \ge r - m)$. First we must bound the expected size of $S$. By Lemma~\ref{chernoff} we have, for each vertex $y \in \Gamma(x)$,
$$\Pr\big(|\Gamma(y) \cap A^{(0)}| \ge r \big) \; \ge \; \exp\left( - 2 \delta^2 n - 4 \delta \sqrt{\frac{n}{\log n}} - \frac{\log\log n}{2} - 6 \right),$$
where $\delta = \ds\frac{1}{2} - p = \frac{1}{2}\sqrt{\frac{\log n}{n}} - \frac{\lambda \log \log n}{\sqrt{n \log n}}$. (Note that we must use Lemma~\ref{chernoff} with $C = 101$ here, since we assume $x \notin A^{(0)}$.) Observe that
$$2 \delta^2 n \, + \, 4 \delta \sqrt{\frac{n}{\log n}} \; = \; \frac{\log n}{2} \, - \, 2\lambda \log\log n \, + \, 2 \, + \, O\left( \frac{(\log \log n)^2}{\log n} \right),$$ and note that, by Lemma~\ref{silly}, $\Pr\big( y \notin A^{(0)} \big| |R| \ge r - m  \big) \ge 1/3$ if $n$ is sufficiently large. Thus
\begin{eqnarray*}
\Ex\left( |S| \, \big| \, |R| \ge r - m \right) & = & n \, \Pr\left( y \notin A^{(0)} \, \big| \, |R| \ge r - m  \right) \Pr\big(|\Gamma(y) \cap A^{(0)}| \ge r \big)\\[+0.5ex]
& \ge & \frac{n}{3e^6} \exp\left( - 2 \delta^2 n - 4 \delta \sqrt{\frac{n}{\log n}} - \frac{\log \log n}{2} \right)\\
& \ge & \frac{ \sqrt{n} \, ( \log n )^{(4\lambda - 1)/2}}{3e^9}
\end{eqnarray*}
since $n$ is sufficiently large, so we can assume the term $O\left( \frac{(\log \log n)^2}{\log n} \right)$ is at most $1$.

Now, we need to show that the variance is not too big. Consider two distinct vertices $y,z \in \Gamma(x)$, and note that $|\Gamma(y) \cap \Gamma(z)| = 2$, and that $\Pr(y \in S \,|\, z \in S) \ge \Pr(y \in S)$. Let $\Gamma(y) \cap \Gamma(z) = \{w,x\}$. Let $p' = \Pr(y \in A^{(0)} \,|\, |R| \ge r - m)$, and note that, by Lemma~\ref{silly}, $p' = \big( 1+o(1) \big) p$. Then
\begin{align*}
& \Pr\big( y \in S \,|\, z \in S, |R| \ge r - m \big) \; \le \; \Pr\big(y \in S \,|\, w \in A^{(0)}, |R| \ge r - m \big)\\[+1ex]
& \hspace{1.5cm} = \; \Pr\big(y \notin A^{(0)} \,|\, |R| \ge r - m\big) \, \Pr\big(|(\Gamma(y) \cap A^{(0)}) \setminus \{x,w\}| \ge r - 1\big)\\[+1ex]
& \hspace{5.37cm} = \; \big(1 - p'\big)\,\Pr\big(S'(n) \ge r\big),
\end{align*}
where $S'(n) \sim 1 + \textup{Bin}(n-2,p)$. But
$$\Pr\big( y \in S \: | \: |R| \ge r - m \big) \; = \; \big( 1 - p' \big) \,\Pr\big( S(n) \ge r \big),$$ where $S(n) \sim \textup{Bin}(n-1,p)$. Thus, by Lemma~\ref{giveone},
$$\Pr(y \in S \,|\, z \in S) \; = \; \big( 1 + o(1) \big)\Pr(y \in S).$$ Finally note that $\Ex(|S|) \to \infty$, and so $\textup{Var}(|S|) = o\big(\Ex(|S|)^2\big)$.

By Lemma~\ref{2mm}, it follows that for any $\eps > 0$, $|S| \ge (1 - \eps)\Ex(|S|)$ with high probability. In particular,
$$\Pr\left(|S| \le \frac{\sqrt{n} \, (\log n)^{(4\lambda - 1)/2}}{e^{11}} \right) \; \le \; \Pr\left(|S| \le \frac{\Ex(|S|)}{2} \right) \; = \; o(1),$$ as claimed.
\end{proof}

Now, let $m = \ds\frac{\sqrt{n} (\log n)^{(4\lambda - 1)/2}}{e^{11}}$, let $\widehat{R}$ denote the event that $|R| \ge r - m$, and let $\widehat{S}$ denote the event that $|S| \ge m$. Since $\lambda > 1/2$, we have $m \gg \sqrt{n \log n}$, and so $r - m < np = \Ex(|R|)$. Therefore, by Lemma~\ref{morethanhalf},
$$\Pr(\widehat{R}) \; \ge \; 1/2 + o(1)$$ as $n \to \infty$. Moreover, by the claim we have
$$\Pr(\widehat{S} \, | \, (x \notin A^{(0)}) \wedge \widehat{R}) \; = \; 1 - o(1).$$ Thus
\begin{eqnarray*}
\Pr(x \in A_r^{(2)}) & \ge & \Pr(x \in A^{(0)}) \: + \: \Pr \big(\widehat{R} \wedge \widehat{S} \, | \, x \notin A^{(0)}) \, \Pr(x \notin A^{(0)})\\[+0.5ex]
& = & p \: + \: (1-p)\Pr \big(\widehat{R} \big) \, \Pr \big(\widehat{S} \, | \, (x \notin A^{(0)}) \wedge \widehat{R})\\[+0.5ex]
& \ge & \frac{3}{4} + o(1),
\end{eqnarray*}
as required.
\end{proof}

\begin{rmk}
Note that we only needed $\lambda > 1/2$ in order to show that
$$m \; = \; \ds\frac{\sqrt{n} (\log n)^{(4\lambda - 1)/2}}{e^{11}} \; \ge \; \sqrt{n \log n} \; \ge \; \left(\ds\frac{1}{2} - p\right) n + 100.$$
Thus, our proof will actually give
$$p_c(Q_n,n/2) \; \le \; \frac{1}{2} \: - \: \frac{1}{2}\sqrt{\frac{\log n}{n}} \: + \: \frac{\log \log n}{2\sqrt{n \log n}} \: + \: O\left(\frac{1}{\sqrt{n \log n}} \right).$$
\end{rmk}

On first sight, it would seem that we are basically done, since in the third round essentially all of the remaining healthy vertices should be infected. Unfortunately, and crucially, however, we have lost independence. The next step, in which we go from most of the vertices to almost all of them, in fact turns out to be the most problematic.

Recall that the elements of $A$ are always chosen independently at random with probability $p$.

\begin{lemma}\label{2to5}
For each $\delta > 0$ there exists a constant $c = c(\delta) > 0$ such that the following holds. Let $n \in \N$, and let $p \in \left( \ds\frac{1}{2} - \sqrt{\ds\frac{\log n }{ n }} \, , \, \ds\frac{1}{2} \right)$ be large enough so that,
$$\Pr\left(x \in A_{n/2+3}^{(2)}\right) \; \ge \; \frac{1}{2} + \delta.$$ Then $\Pr(x \in A^{(5)}) \ge 1 - e^{-c n}$.
\end{lemma}

\begin{proof}
We wish to estimate $\Pr(x \notin A^{(5)})$; in order to do so, we must show that the probability that $|\Gamma(x) \setminus A^{(4)}| \ge n/2$ is small. We begin with an important fact:\\[-1.5ex]

\hspace{0.1cm} If $S \subset Q_n$ satisfies $d(y,z) \ge 5$ for every $y,z \in S$ with $y \neq z$,

\hspace{0.1cm} then, for any $r$, the events $\{y \in A_r^{(2)}\}_{y \in S}$ are independent.\\[-1.5ex]

\noindent Let $Y = S(x,3) = \{y \in Q_n : d(x,y) = 3\}$, so $|Y| = {n \choose 3}$, and let $m = 3{n \choose 2} \le 2n^2$. By Lemma~\ref{hyperpartition} there exist disjoint sets $B_1, \ldots, B_m$ such that $\bigcup B_j = Y$, and for each $j \in [m]$, the events $\{y \in A^{(2)}\}_{y \in B_j}$ are independent.

Now, let $\eps > 0$, and with foresight, observe that by
Lemma~\ref{normalcher},
\begin{equation}\label{foresight}
\Pr\left(|\Gamma(x) \cap A^{(0)}| \,<\, \frac{n}{2} \,-\, \delta^2\eps^2 n \right) \; \le \; \exp\left( - \delta^4\eps^4 n \right),
\end{equation}
since $\left(\ds\frac{1}{2} - p\right)n \ll \delta^2\eps^2 n$. Let $J(x)$ denote the event that this does not happen, i.e., that $|\Gamma(x) \cap A^{(0)}| \ge n/2 - \delta^2\eps^2 n$, and assume from now on that $J(x)$ holds. Now choose $S \subset \Gamma(x) \setminus A^{(0)}$ with $|S| = n/2$. We thus have at most ${{n/2 + \delta^2\eps^2 n} \choose n/2}$ choices for $S$. We shall show that it is extremely unlikely that $S \cap A^{(4)} = \emptyset$.

Indeed, let $T = \Gamma(S) \setminus \{x\}$, so $d(x,y) = 2$ for every $y \in T$, and note that $T = T_1 \cup T_2$, where $T_i = \{y \in T : |\Gamma(y) \cap S| = i\}$ for $i = 1,2$. Note further that $|T_1| = n^2/4$ and $|T_2| = {{n/2} \choose 2}$.

Let $a = |T_1 \cap A^{(3)}|$ and $b = |T_2 \cap A^{(3)}|$, and suppose that $S \cap A^{(4)} = \emptyset$. Then $a + 2b \le n^2/4$. (This follows by simply counting edges -- each vertex of $S$ sends at least $n/2$ edges to vertices not in $A^{(3)}$.) But, by Lemma~\ref{normalcher}, since $p \in \left( \ds\frac{1}{2} - \sqrt{\ds\frac{\log n }{ n }} \, , \, \ds\frac{1}{2} \right)$, we have
\begin{equation}\label{T1}
\Pr\left( \left| |T_1 \cap A^{(0)}| \,-\, \frac{n^2}{8} \right| \, \ge \, 2\sqrt{n^3 \log n} \right) \: \le \: \exp\left( - n \log n \right),
\end{equation}
and similarly
\begin{equation}\label{T2}
\Pr\left( \left| |T_2 \cap A^{(0)}| \,-\, \frac{n^2}{16} \right| \, \ge \, 2\sqrt{n^3 \log n} \right) \: \le \: \exp\left( - n \log n \right).
\end{equation}
Let $F(S)$ denote the event that neither of these events occurs, i.e., that $\left| |T_i \cap A^{(0)}| \,-\, n^2/8i \right| \le 2\sqrt{n^3 \log n}$ for $i = 1,2$, so $\Pr(F(S)) \ge 1 - 2n^{-n}$, and assume that the event $F(S)$ holds for every $S$.

So, with very high probability, if $|T \cap A^{(3)} \setminus A^{(0)}| > 6\sqrt{n^3 \log n}$ then $a + 2b > n^2/4$, and thus $S \cap A^{(4)} \neq \emptyset$. Hence the following claim will (essentially) complete the proof.\\[-1ex]

\noindent \ul{Claim}:
\begin{equation}\Pr\left( |T \cap A^{(3)} \setminus A^{(0)}| = O\left( \sqrt{n^3 \log n} \right) \right) \: \le \: 6n^2 \exp\left( - \frac{\delta^2 \eps n}{2} \right).\label{claim}\end{equation}

\begin{proof}[Proof of claim]
Consider the bipartite graph $H$, with vertex set $W \cup Y$, where $W = T \setminus A^{(0)}$ and $Y = S(x,3)$, and edge set $\{wy : wy \in E(Q_n)\}$. Furthermore, let us colour the edges of $H$ red and blue, according to whether or not the endpoint in $Y$ is also in $A^{(2)}$, i.e., $c(wy) =$ red if and only if $y \in A^{(2)}$.

Now, $e(H) = (n-2)|W|$, and $|W| \, =\, \ds\frac{3n^2}{16} \, + \, O(\sqrt{n^3\log n})$, since we assumed that $F(S)$ holds. Suppose $|T \cap A^{(3)} \setminus A^{(0)}| = O\left( \sqrt{n^3 \log n} \right)$. Then only $O\left( \sqrt{n^3 \log n} \right)$ vertices of $W$ have at least $n/2$ neighbours in $A^{(2)}$, and thus the number of red edges in $H$, $e_R(H)$, satisfies
$$e_R(H) \; \le \; \frac{n|W|}{2} \,+\, O(\sqrt{n^5 \log n}) \; = \; \frac{3n^3}{32} \,+\, O(\sqrt{n^5 \log n}),$$
and thus
$$e_B(H) \; \ge \; \frac{n|W|}{2} \,-\, O(\sqrt{n^5 \log n}) \; = \; \frac{3n^3}{32} \,-\, O(\sqrt{n^5 \log n})$$
also. Now, recall the partition $B_1, \ldots, B_m$ of $Y$ into independent sets, and define a refinement of it by setting
$$B_j(i) \; = \; \{z \in B_j \,:\, |\Gamma(z) \cap T| = i\}$$ for each $i \in [3]$ and each $j \in [m]$. Note that the sets $B_j(i)$ are determined by the set $S$ and the partition $B_1, \ldots, B_m$.\\[-1ex]

\noindent \ul{Subclaim}: If $\eps > 0$ (chosen above) is small enough, then there exists an $i \in [3]$ and a $j \in [m]$ such that\\[-1.5ex]

$(a)$ $\ell \ge 3\eps n$ edges of $H$ are incident with $B_j(i)$, and\\[-2ex]

$(b)$ at most $(1/2 + \delta/2)\ell$ of these edges are red.

\begin{proof}[Proof of subclaim]
Suppose the subclaim is false. Let us count the total number of blue edges. From those $B_j(i)$ with at most $3\eps n$ incident edges, we get at most $9\eps m n \le 18\eps n^3$ edges. From the others, we get at most $e(H)(1/2 - \delta/2) \le n|W|(1/2 - \delta/2)$ blue edges. If $\eps$ is small enough, this contradicts the bound on $e_B(H)$ above.
\end{proof}

Now, recall that the events $\{y \in A^{(2)}\}_{y \in B_j(i)}$ are independent, and that we chose $p$ so that $\Pr(y \in A_{n/2+3}^{(2)}) \ge 1/2 + \delta$ for each $y \in Q_n$. Unfortunately, this event is not independent of which members of $B(x,2)$ are in $A^{(0)}$; however, for each vertex $y \in B_j(i)$, $y$ has at most three neighbours in $B(x,2)$. Considering the $(n-3)$-dimensional sub-hypercube containing $y$ but none of these neighbours, we see that, for any set $B(x,2) \cap A^{(0)}$, we have $\Pr(y \in A_{n/2}^{(2)}) \ge 1/2 + \delta$.

Let $E_j(i)$ denote the event that the set $B_j(i)$ satisfies the conditions $(a)$ and $(b)$ of the subclaim, and note again that each vertex in $B_j(i)$ sends $i \le 3$ red edges into $T$. Thus, by Lemma~\ref{normalcher},
\begin{eqnarray}
\Pr\left( \bigcup_{i,j} E_j(i) \right) & \le & \sum_{i,j} \Pr\left( \textup{Bin}\left(\frac{3\eps n}{i},\frac{1}{2} + \delta \right) \, \le \, \left( \nonumber \frac{1}{2} + \frac{\delta}{2} \right) \frac{3\eps n}{i} \right) \\[+1ex]
& \le & 6n^2 \exp\left( - \frac{\delta^2 \eps n}{2} \right).\nonumber
\end{eqnarray}
However, we showed that if $|T \cap A^{(3)} \setminus A^{(0)}| = O\left( \sqrt{n^3 \log n} \right)$, then one of the events $E_j(i)$ occurs. Hence this proves the claim.
\end{proof}

Let $M(S)$ denote the event that $|T \cap A^{(3)} \setminus A^{(0)}| > 6\sqrt{n^3 \log n}$, so by the claim, $\Pr(M(S)) \,\ge\, 1 - 6n^2 \exp\left( - \ds\frac{\delta^2 \eps n}{2} \right)$. Finally, recall that, assuming the event $J(x)$ holds, we had at most
\begin{equation}\label{choices}
\ds{{n/2 + \delta^2\eps^2n} \choose \delta^2\eps^2 n} \: \le \: \left( \frac{e}{\delta^2\eps^2} \right)^{\delta^2\eps^2 n} \: \le \: \exp\left( \frac{\delta^2 \eps n}{10} \right)
\end{equation}
choices for the set $S$ if $\eps$ is sufficiently small, since $(\frac{e}{x^2} )^x < e^{\delta/10}$ if $x$ is sufficiently small.

Now, we claim that if each of the events $J(x)$, $F(S)$ and $M(S)$ holds (for all $S$ as above), then $x \in A^{(5)}$. Indeed, $F(S) \wedge M(S)$  implies that $a + 2b > n^2/4$, and therefore that $S \cap A^{(4)} \neq \emptyset$, as explained above. If this holds for every $S$, then it follows that $x \in A^{(5)}$. Therefore, by (\ref{foresight}), (\ref{T1}), (\ref{T2}), (\ref{claim}) and (\ref{choices}),
\begin{eqnarray*}
\Pr(x \notin A^{(5)}) & \le & \Pr(J(x)^c) \: + \: \sum_S \Pr(F(S)^c) \: + \: \sum_S \Pr(M(S)^c)\\
& \le & \exp\left( - \delta^4\eps^4 n \right) \; + \; \exp\left( \frac{\delta^2 \eps n}{10} \right) \left( 6n^2 \exp\left( - \frac{\delta^2 \eps n}{2}\right) + 2n^{- n} \right)\\[+1ex]
& \le & \exp\big( - c n \big)
\end{eqnarray*}
for some $c = c(\delta) > 0$, as required.
\end{proof}

Finally, we make the (easier) jump from exponential to super-exponential probability; this step takes us from round $k$, to round $2k+1$.

\begin{lemma}\label{exptoall}
Let $k \in \N$, and let $n$ be sufficiently large (in particular, $n \ge 4k$). Let $x \in Q_n$ and $c > 0$. Suppose that $p$ is chosen such that $$\Pr(y \in A^{(k)}) < e^{-c n}$$ for each $y \in Q_n$. Then there exists a constant $d = d(c,k) > 0$ such that
$$\Pr(x \in A^{(2k+1)}) \; \le \; e^{-dn^2}.$$
\end{lemma}

\begin{proof}
As in Lemma~\ref{2to5}, note first that\\[-1.5ex]

\hspace{0.1cm} if $S \subset Q_n$ satisfies $d(y,z) \ge 2k+1$ for every $y,z \in S$ with $y \neq z$,

\hspace{0.1cm} then the events $\{y \in A^{(k)}\}_{y \in S}$ are independent.\\[-1.5ex]

\noindent Let $x \in Q_n$, and let $m = (k+1){n \choose k} \le 2n^k$. By Lemma~\ref{hyperpartition}, there exist disjoint sets $B_1, \ldots, B_m$ such that $\bigcup B_j = S(x,k+1)$, and for each $j \in [m]$, the events $\{y \in A^{(k)}\}_{y \in B_j}$ are independent.

The argument is now very simple. Suppose $x \notin A^{(2k+1)}$; we claim that for each $t \in [0,k+1]$ there exist a set $T(t) \subset S(x,t)$ such that $T(t) \cap A^{(2k+1-t)} = \emptyset$, and
$$|T(t)| \; \ge \; \frac{(n/2)!}{(n/2 - t)! t!} \; \ge \; \frac{n^t}{4^t t!}.$$

Indeed, let $T(0) = \{x\}$, and note that since $x \notin A^{(2k+1)}$, $T(0)$ satisfies the conditions. Now, suppose we have found $T(t)$ as required. Then each $y \in T(t)$ has at most $n/2$ neighbours in $S(x,t+1) \cap A^{(2k-t)}$, and so at least $n/2 - t$ neighbours in $S(x,t+1) \setminus A^{(2k-t)}$. Moreover, each element of $S(x,t+1)$ has exactly $t+1$ neighbours in $S(x,t)$. Thus, by counting edges, there must exist a set $T(t+1) \subset S(x,t+1)$ with $T(t+1) \cap A^{(2k-t)} = \emptyset$, and $|T(t+1)| \ge (n/2 - t)|T(t)|/(t+1)$, as required. The second inequality holds since $n \ge 4k$.

Consider $T(k+1)$, and note that it has at least $\alpha n^{k+1}$ elements, where $\alpha > 0$ does not depend on $n$. Thus, there must exist a $j \in [m]$ and an absolute constant $\eps > 0$ (not depending on $n$) such that $|B_j| \ge \eps n$, and $|T(k+1) \cap B_j| \ge \eps |B_j|$. Indeed, if no such $j$ exists, then we would have $|T(k+1)| \le \eps mn + \eps {n \choose k+1} < \alpha n^{k+1}$ if $\eps$ is sufficiently small (recall that $m \le 2n^k$).

Now, recall that $\Pr(y \notin A^{(k)}) < e^{-c n}$ for each $y \in B_j$, that the events $\{y \in A^{(k)}\}_{y \in B_j}$ are independent, and that $T(k+1) \cap A^{(k)} = \emptyset$. Thus, by Lemma~\ref{nunlikely} we have
$$\Pr\left(|T(k+1) \cap B_j| \ge \eps |B_j| \, \big| \, |B_j| \ge \eps n \right) \; \le \; e^{-b|B_j|^2} \; \le \; e^{-b\eps^2n^2}$$ for some $b = b(c,\eps) > 0$, and so
\begin{eqnarray*}
\Pr(x \notin A^{(2k+1)}) & \le & \Pr(\exists j\textup{ with } |B_j| \ge \eps n\textup{ and }|T(k+1) \cap B_j| \ge \eps |B_j|)\\
& \le & me^{-b\eps^2n^2} \; \le \; e^{-dn^2}
\end{eqnarray*} for some $d = d(c, \eps) > 0$. Since $\eps = \eps(k)$, the lemma follows.
\end{proof}

Before moving on to the lower bound, let us put together the pieces from this section.

\begin{cor}\label{upperbound}
Let $\lambda > 1/2$, let $n \in \N$ be sufficiently large, and let $$p(n) \; = \; \frac{1}{2} \: - \: \frac{1}{2}\sqrt{\frac{\log n}{n}} \: + \: \frac{\lambda \log \log n}{\sqrt{n \log n}}.$$ Then, in majority bootstrap percolation on $Q_n$, with initial set $A$ of density $p$,
$$\Pr(A\textup{ percolates}) \; \to \; 1$$
as $n \to \infty$.
\end{cor}

\begin{proof}
Let $n$ and $p$ be as given. By Lemma~\ref{3/4} it follows that $\Pr(x \in A^{(2)}) \ge 2/3$, for sufficiently large $n$ (depending on $\lambda$). So, by Lemma~\ref{2to5}, it follows that $\Pr(x \notin A^{(5)}) \le e^{-cn}$ for some $c = c(\lambda) > 0$, and thus by Lemma~\ref{exptoall}, it follows that $\Pr(x \notin A^{(11)}) \le e^{-dn^2}$ for some $d = d(\lambda) > 0$. Finally, the union bound gives
\begin{eqnarray*}
\Pr(A\textup{ does not percolate}) & \le & \Pr\left( \bigcup_{x \in Q_n} \big(x \notin A^{(11)}\big)\right) \\
& \le & \sum_{x \in Q_n} \Pr \big(x \notin A^{(11)}\big) \; \le \; 2^n e^{-dn^2} \; = \; o(1),
\end{eqnarray*} as required.
\end{proof}

\section{Proof of the lower bound in Theorem~\ref{sharp}}\label{lowersec}

Let $n \in \N$ be sufficiently large, and let $$p(n) \; = \; \frac{1}{2} \: - \: \frac{1}{2}\sqrt{\frac{\log n}{n}} \: + \: \frac{\lambda \log \log n}{\sqrt{n \log n}},$$ for some $\lambda \in \RR$, as in the previous section. We shall couple the bootstrap process with a modified process in which, if $\lambda \le -2$, then $A^{(4)} = A^{(3)} \neq Q_n$ with high probability.

We shall start slowly, and build up to the full result. First, we consider just a two step process, and show that if $p \le 1/2 - \eps$ then the original process does not percolate. Let us refer to the original (majority) process as \textbf{Boot}, and let $t = t(n) \ge 0$ be any non-negative function. We define the process \textbf{Boot1($t$)} as follows.\\[-1ex]
\begin{itemize}
\item The elements of $A^{(0)}$ are chosen independently at random, each with probability $p$.\\[-1.5ex]
\item $x \in A^{(1)}$ if $x \in A^{(0)}$ or $|\Gamma(x) \cap A^{(0)}| \ge n/2 - t$.\\[-1.5ex]
\item If $i \ge 1$, then $x \in A^{(i+1)}$ if $x \in A^{(i)}$ or $|\Gamma(x) \cap A^{(i)}| \ge n/2$.\\[-1ex]
\end{itemize}
Note that \textbf{Boot} = \textbf{Boot($0$)}, and that the process \textbf{Boot1($t$)} dominates the process \textbf{Boot}, in the sense that given the same initial set $A^{(0)}$, then for each $k \in \N$, the set $A^{(k)}$ given by \textbf{Boot1($t$)} contains that given by \textbf{Boot}.

The following simple result, together with Corollary~\ref{upperbound}, implies that the critical probability for percolation in the hypercube is $1/2 + o(1)$.

\begin{prop}\label{boot1}
Let $\eps > 0$, and suppose $p = 1/2 - \eps$ and $t = \eps n/4$. Then in \textup{\textbf{Boot1($t$)}}, $A^{(2)} = A^{(1)}$ with high probability.
\end{prop}

\begin{proof}
Let $x \in Q_n$, and suppose $x \in A^{(2)} \setminus A^{(1)}$. Then there must exist a set $T \subset \Gamma(x)$, with $|T| = t$, and $T \subset A^{(1)} \setminus A^{(0)}$. This is because $x \in A^{(2)} \setminus A^{(1)}$ implies $|\Gamma(x) \cap A^{(1)}| \ge n/2$, and $x \notin A^{(1)}$ implies $|\Gamma(x) \cap A^{(0)}| < n/2 - t$. We shall show that $\Pr(T$ exists$) < e^{-cn^2}$ for some $c > 0$.

Indeed, recall that $S(x,k) = \{y \in Q_n : d(x,y) = k\}$ for each $k \in \N$, and consider the set $\Gamma(T) \cap S(x,2)$. It has ${t \choose 2}$ elements with two neighbours in $T$, and $t(n-t)$ elements with one neighbour in $T$. Denote these two sets $B$ and $C$ respectively. Now, we claim that since $T \subset A^{(1)} \setminus A^{(0)}$, we have
\begin{equation}|C \cap A^{(0)}| \; \ge \; t\left(\frac{n}{2} - t\right) \: - \: 2{t \choose 2} \; \ge \; \frac{nt}{2} \: - \: 2t^2.\label{CA0}\end{equation}
This follows by counting edges. Indeed, note that each member of $T$ has at least $n/2 - t$ neighbours in $A^{(0)}$, and we assumed that $x \notin A^{(0)}$. Each vertex of $C$ has only one neighbour in $T$, so even if every member of $B$ is in $A^{(0)}$ we still get the bound (\ref{CA0}). But \begin{eqnarray*}
\Ex\left( |C \cap A^{(0)}| \right) & = & pt(n-t) \; = \; \frac{nt}{2} \: - \: \eps nt \: - \: \left(\frac{1}{2} \, - \, \eps\right)t^2\\
& \le & \frac{nt}{2} \, - \, \eps n t,
\end{eqnarray*}
and $|C \cap A^{(0)}| \sim \textup{Bin}(t(n-t),p)$, so, since $t = \eps n/4$,
\begin{eqnarray*}
\Pr\left(|C \cap A^{(0)}| \ge \frac{nt}{2} - 2t^2 \right) & \le & \exp\left( -\frac{2(\eps n t - 2t^2)^2}{t(n-t)} \right) \; \le \;
\exp\left( -\frac{\eps^3 n^2}{8} \right)
\end{eqnarray*}
by Lemma~\ref{normalcher}. But we have only at most $2^{n}$ choices for the set $T$, so
$$\Pr\big( x \in A^{(2)} \setminus A^{(1)} \big) \; \le \;  2^n \exp\left( -\frac{\eps^3 n^2}{8} \right) \; < \; e^{-cn^2}$$
for some $c > 0$, and so, since $|V(Q_n)| = 2^n$,
$$\Pr\big(A^{(2)} \setminus A^{(1)} \neq \emptyset \big) \; \le \; 2^n \, \Pr\big( x \in A^{(2)} \setminus A^{(1)} \big) \; = \; o(1),$$ as required.
\end{proof}

The following corollary is immediate from the proposition and Corollary~\ref{upperbound}.

\begin{cor}
$p_c(Q_n,n/2) = \ds\frac{1}{2} + o(1)$.
\end{cor}

\begin{proof}
By Corollary~\ref{upperbound} we have $p_c(Q_n,n/2) \le 1/2 + o(1)$, so let $\eps > 0$, $t = \eps n/4$ and $p = 1/2 - \eps$, and consider the \textbf{Boot1($t$)} process on $Q_n$. First note that $A^{(1)} \neq V(Q_n)$ with high probability, since, by Lemma~\ref{normalcher}, $\Pr\big(|A^{(1)} \setminus A^{(0)}| \ge 2^n/100 \big) = o(1)$. Thus, by Proposition~\ref{boot1}, the process \textbf{Boot1($t$)} does not percolate with high probability.

Now, coupling \textbf{Boot} with \textbf{Boot1($t$)} in the obvious way, we see that also \textbf{Boot} does not percolate whp, so $p_c(Q_n,n/2) \ge 1/2 - \eps$.
\end{proof}

\begin{rmk}
In fact, letting $\eps = \alpha n^{-1/3}$ for some $\alpha \in \RR$ with $\alpha^3 > 32$, the same proof in fact gives
$$p_c(Q_n,n/2) \; > \; \frac{1}{2} \,-\, \frac{\alpha}{n^{1/3}},$$
since $\Pr(x \in A^{(2)} \setminus A^{(1)}) \le 2^n e^{-\alpha^3n/8}$, and so $\Pr(A^{(2)} \setminus A^{(1)} \neq \emptyset) = o(1)$.
\end{rmk}

In order to improve this bound, we have to allow the process to go a little further. We call the following process \textbf{Boot3($t$)}.\\[-1ex]
\begin{itemize}
\item The elements of $A^{(0)}$ are chosen independently at random.\\[-1.5ex]
\item $x \in A^{(1)}$ if $x \in A^{(0)}$ or $|\Gamma(x) \cap A^{(0)}| \ge n/2 - 3t$.\\[-1.5ex]
\item $x \in A^{(2)}$ if $x \in A^{(1)}$ or $|\Gamma(x) \cap A^{(1)}| \ge n/2 - 2t$.\\[-1.5ex]
\item $x \in A^{(3)}$ if $x \in A^{(2)}$ or $|\Gamma(x) \cap A^{(2)}| \ge n/2 - t$.\\[-1.5ex]
\item If $i \ge 3$, then $x \in A^{(i+1)}$ if $x \in A^{(i)}$ or $|\Gamma(x) \cap A^{(i)}| \ge n/2$.\\[-1ex]
\end{itemize}

Note that, since we are trying to distinguish between values of $p$ which differ by $O\left(\ds\frac{\log \log n}{\sqrt{n \log n}}\right)$, we should take $t$ no larger than $n$ times this. In fact, we shall show that if $t = \ds\sqrt{\frac{n}{\log n}}$, $p$ is as above and $\lambda \le -2$, then in the \textbf{Boot3($t$)} process we have $A^{(4)} = A^{(3)} \neq Q_n$ with high probability. We begin by showing that $A^{(3)} \neq Q_n$ (in fact we only need the slightly weaker condition, that $\lambda < 1/4$).

\begin{lemma}\label{notall}
Let $\lambda < 1/4$, let $\delta = \ds\frac{1}{2} \sqrt{\ds\frac{\log n}{n}} - \frac{\lambda \log \log n}{\sqrt{n \log n}}$, and suppose that $p = \ds\frac{1}{2} - \delta$ and $t = \ds\sqrt{\frac{n}{\log n}}$. Then in \textup{\textbf{Boot3($t$)}}, $A^{(3)} \neq Q_n$ whp.
\end{lemma}

\begin{proof}
Let $x \in Q_n$. We shall show that $\Pr(x \notin A^{(3)}) > e^{-cn}$ for some small $c > 0$, and that there is a set $\{x_1,\ldots,x_\ell\}$, with $\ell \ge 2^n / n^6$, for which the statements $x_i \in Q_n$ are independent. The proof that $\Pr(x \in A^{(3)})$ is not too big is similar to that of Lemma~\ref{3/4}.

First we show that $\Pr(y \in A^{(2)} \, | \, y \notin A^{(0)}) = o(1)$, for each $y \in \Gamma(x)$. In fact, we shall need a slightly stronger result: that this still holds, given any set $\Gamma(x) \cap A^{(0)}$. Clearly the events $y \in A^{(2)}$ and $z \in A^{(0)}$ are positively correlated, so let $y \in \Gamma(x) \setminus A^{(0)}$, and assume that $\Gamma(x) \cap A^{(0)} = \Gamma(x) \setminus \{y\}$. Now, let $R = \Gamma(y) \cap A^{(0)}$ and $S = \Gamma(y) \cap A^{(1)} \setminus A^{(0)}$. We want to show that $\Pr(|R| + |S| \ge n/2 - 2t) \to 0$ as $n \to \infty$. This follows easily from the following claim.\\[-1ex]

\noindent\ul{Claim}: $\Pr\left( |S| > e^{10}\sqrt{n} \,(\log n)^{2\lambda} \right) = o(1).$

\begin{proof}[Proof of claim]
We use the second moment method (Lemma~\ref{2mm}). First we must bound the expected size of $S$. By Lemma~\ref{normalcher} we have, for each vertex $z \in \Gamma(y)$,
\begin{eqnarray*}
\Pr\left(|\Gamma(z) \cap A^{(0)}| \ge \frac{n}{2} - 3t \right) & = & \Pr\left(|\Gamma(z) \cap A^{(0)} \setminus \Gamma(x)| \ge \frac{n}{2} - 3t - 1 \right)\\[+0.5ex]
& \le & \exp\left( -\frac{2 (\delta n - 4t)^2 }{n} \right),
\end{eqnarray*}
where $\delta = \ds\frac{1}{2} - p = \frac{1}{2}\sqrt{\frac{\log n}{n}} - \frac{\lambda \log \log n}{\sqrt{n \log n}}$, since $z$ has only one other neighbour in $\Gamma(x)$, and $t \ge 1$. Observe that
$$2 \delta^2 n \, - \, 16 \delta t + \frac{32t^2}{n} \; = \; \frac{\log n}{2} \, - \, 2\lambda \log\log n \, - \, 8 \, + \, O\left( \frac{(\log \log n)^2}{\log n} \right),$$ and so
\begin{eqnarray*}
\Ex(|S|) & = & n \, \Pr\big( z \notin A^{(0)} \big) \Pr\big(|\Gamma(z) \cap A^{(0)}| \ge n/2 - 3t\big)\\[+1ex]
& \le & n \exp\left( -\frac{\log n}{2} \, + \, 2\lambda \log\log n \, + \, 9 \right) \; = \; e^9 \sqrt{n} \, ( \log n )^{2\lambda}
\end{eqnarray*}
since $n$ is sufficiently large, so we can assume the term $O\left( \frac{(\log \log n)^2}{\log n} \right)$ is at most $1$. It is similarly straightforward (using Lemma~\ref{chernoff}) to show that $\Ex(|S|) \to \infty$ as $n \to \infty$.

Now, we need to show that the variance is not too big. Consider two distinct vertices $u,v \in \Gamma(y) \setminus \{x\}$, and note that $|\Gamma(u) \cap \Gamma(v)| = 2$, and that $\Pr(u \in S \,|\, v \in S) \ge \Pr(u \in S)$. Let $\Gamma(u) \cap \Gamma(v) = \{w,y\}$, and also let $\Gamma(u) \cap \Gamma(x) = \{a,y\}$. Then, since we assumed $y \notin A^{(0)}$,
$$\Pr(u \in S \,|\, v \in S) \; \le \; \Pr(u \in S \,|\, w \in A^{(0)}) \; = \; (1-p)\,\Pr(S'(n) \ge n/2-3t),$$
where $S'(n) = |\left( \Gamma(u) \cap A^{(0)} \right) \cup \{a,w\} \setminus \{y\}| \sim 2 + \textup{Bin}(n-3,p)$. But
$$\Pr(u \in S) = (1-p)\,\Pr(S(n) \ge n/2-3t),$$ where $S(n) = |\left( \Gamma(u) \cap A^{(0)} \right) \cup \{a\} \setminus \{y\}| \sim 1 + \textup{Bin}(n-2,p)$. Thus, by Lemma~\ref{giveone},
$$\Pr(u \in S \,|\, v \in S) \; = \; \big( 1 + o(1) \big)\Pr(u \in S),$$ and so, since $\Ex(|S|) \to \infty$ as $n \to \infty$, we have $\textup{Var}(|S|) = o\big(\Ex(|S|)^2\big)$.

By Lemma~\ref{2mm}, it follows that $|S| \le 2\Ex(|S|)$ with high probability. Thus
$$\Pr\left( |S| > e^{10}\sqrt{n} \,(\log n)^{2\lambda} \right) \; \le \; \Pr\left(|S| > 2\Ex(|S|) \right) \; = \; o(1),$$ as claimed.
\end{proof}

Now, let $m = e^{10}\sqrt{n} \,(\log n)^{2\lambda}$, let $\widehat{R}$ denote the event that $|R| \ge n/2 - \delta n / 2$, and let $\widehat{S}$ denote the event that $|S| \ge \delta n /2 - 2t \ge \delta n /3$. Note that $\delta n \gg m$ since $\lambda < 1/4$. Thus $\Pr(\widehat{R}) = o(1)$, by Lemma~\ref{normalcher}, and $\Pr(\widehat{S}) = o(1)$, by the claim, and hence
\begin{eqnarray*}
\Pr(y \in A^{(2)} \, | \, y \notin A^{(0)}) & \le & \Pr(\widehat{R}) \: + \: \Pr \big(\widehat{S}) \; = \; o(1).
\end{eqnarray*}

To complete the proof of the lemma, we need to use the simple fact that for positively correlated events $E_1, \ldots, E_k$,
$$\Pr\left(\bigcap_{i=1}^k E_i \right) \; \ge \; \prod_{i=1}^k \Pr\big(E_i\big).$$ Given $y \in Q_n$, let $T(y)$ denote the event that $y \notin A^{(2)}$. Noting that the events $\{T(y)\}_{y \in \Gamma(x) \setminus A^{(0)}}$ are positively correlated, we have, for any set $\Gamma(x) \cap A^{(0)}$, any (small) $\eps > 0$, and sufficiently large $n$,
$$\Pr\big( \Gamma(x) \cap A^{(2)} \setminus A^{(0)} = \emptyset \big) \; = \; \Pr\left(\bigcap_{y \in \Gamma(x) \setminus A^{(0)}} T(y)\right) \; \ge \; (1 - \eps)^n,$$
and therefore
\begin{eqnarray*}
\Pr(x \notin A^{(3)}) & \ge & \Pr(x \notin A^{(0)}) \Pr\left(|\Gamma(x) \cap A^{(0)}| \le \frac{n}{2} - 3t \right) \Pr\left(\bigcap_{y \in \Gamma(x) \setminus A^{(0)}} T(y)\right) \\
& \ge & \frac{1}{2} \times (1 - \eps) \times (1 - \eps)^n \; \ge \; \exp\left(- 2\eps n \right).
\end{eqnarray*}

Finally, if $d(u,v) \ge 7$ then the events $u \in A^{(3)}$ and $v \in A^{(3)}$ are independent, so, by Lemma~\ref{partition}, there exists a set $K$ of size at least $2^n/n^6$ for which the events $\{x \in A^{(3)}\}_{x \in K}$ are all independent. Thus
$$\Pr(x \in A^{(3)}\textup{ for all }x \in K) \le (1 - e^{-2\eps n})^{|K|} \; = \; o(1)$$ if $\eps$ is chosen to be sufficiently small. This completes the proof.
\end{proof}

Next, we shall show that the \textup{\textbf{Boot3($t$)}} process stops after at most three steps if $\lambda \le -2$. We shall need the following lemma about counting 3-uniform hypergraphs.

Given a 3-uniform hypergraph $H$, and $i,j \in [n]$ with $i < j$, we shall write $d_H(i,j)$ for the degree of the pair $\{i,j\}$ in $H$, i.e., $$d_H(i,j) \; = \; |\{k : ijk \in E(H)\}|.$$ We write $H \le G$ if $H$ is a (not necessarily induced) sub-hypergraph of $G$, and define $\|H\| = \sum_{i<j} {{d_H(i,j)} \choose 2}$.

\begin{lemma}\label{count}
Let $G$ be a labelled $3$-uniform hypergraph with $n$ vertices. Further, let $t^2 \le n \ll t^3$, $s \le t^3$, $m \in \N$, and
$$S(G,m,s) \; = \; \{ H \le G : e(H) = s\textup{ and }\|H\| \ge m \}.$$ Then, for sufficiently large $n$,
$$|S(G,m,s)| \; \le 2^n \left( \frac{20nt}{s} \right)^{s} {{e(G)} \choose {s - 2m/n + 2t^5/n}}.$$
\end{lemma}

\begin{proof}
Let $H \in S(G,m,s)$, so $H$ $\le$ $G$, $e(H) = s$ and $\|H\|$ $\ge$ $m$. We partition the elements $K = \{(i,j) : i < j\}$ into two sets, which we imaginatively title `big' and `little'. To be precise, let $\tau = \tau(t)$ be a function to be determined later, let
$$K_1 = \{(i,j) : d_H(i,j) \ge \tau\},$$ be the big set, and let $K_2 = K \setminus K_1$ be the small one. Let $|K_1| = k_1$ and $|K_2| = k_2$, let
$$L = \{e \in E(H) : (i,j) \in K_1\textup{ for some }i,j \in e\},$$ and $|L| = \ell$. We shall show that if $\|H\|$ is large, then $\ell$ must be large, and hence that there are only few choices for $H$.

Indeed, first note that $\sum_{i<j} d_H(i,j) = 3s$, since each edge contains exactly three elements of $K$, and hence $k_1 \le \ds\frac{3s}{\tau}$. Also, note that $\ds\sum_{(i,j) \in K_1} d_H(i,j) \le \ell + {k_1 \choose 2}$, since $\ell$ edges intersect $K_1$, and each pair of elements of $K_1$ is contained in at most one edge $e \in E(H)$. Thus we gain at most ${k_1 \choose 2}$ extra in the sum from those edges of $L$ which contain more than one element of $K_1$.

Now we reach the crux. Recall that $d_H(i,j) \le n$ for all $i$ and $j$, and $d_H(i,j) \le \tau$ if $(i,j) \in K_2$. Thus, by the convexity of ${x \choose 2}$ (Observation~\ref{convexity}), and letting $\tau = t^2/2$, we have
\begin{eqnarray}
\|H\| & \le & \left( \frac{1}{n} \sum_{(i,j) \in K_1} d_H(i,j) \right) {n \choose 2} \: + \: \left( \frac{1}{\tau} \sum_{(i,j) \in K_2} d_H(i,j) \right) {\tau \choose 2}\nonumber\\
& \le & \frac{(\ell + k_1^2)n}{2} \: + \: \frac{3s\tau}{2} \; \le \; \frac{n \ell}{2} \, + \, t^5\label{ell}
\end{eqnarray} since $k_1 \le \ds\frac{3s}{\tau} \le 6t$, so $k_1^2n \ll t^5$.

Now we have only to count. To determine $H$, it suffices to choose the $k_1$ pairs $(i,j)$ in $K_1$, and $\ell$ edges incident to these pairs, and then $s - \ell$ other edges. Thus,
\begin{eqnarray*}
|S(G,m,s)| \; \le \; \sum_{k_1,\ell} {{n^2} \choose k_1} {{k_1n} \choose \ell} {{e(G)} \choose {s - \ell}} \; \le \; 6t \sum_{\ell} n^{12t} \left( \frac{6ent}{s} \right)^{s} {{e(G)} \choose {s - \ell}}
\end{eqnarray*}
since $k_1 \le 6t$, and using the trivial bound $\ell \le s$. But $\ds{{e(G)} \choose {s - \ell}}$ is decreasing in $\ell$, and $\ell \ge \ds\frac{2(m - t^5)}{n}$ by (\ref{ell}), so
\begin{eqnarray*}
|S(G,m,s)| & \le & 6stn^{12t} \left( \frac{6ent}{s} \right)^{s} {{e(G)} \choose {s - 2(m - t^5)/n}}\\[+0.5ex]
& \le & 2^n \left( \frac{20nt}{s} \right)^{s} {{e(G)} \choose {s - 2m/n + 2t^5/n}}
\end{eqnarray*}
for sufficiently large $n$, as claimed.
\end{proof}

We are ready to prove the key lemma.

\begin{lemma}\label{boot3}
Let $\lambda \le -2$, let $\delta = \ds\frac{1}{2} \sqrt{\ds\frac{\log n}{n}} - \frac{\lambda \log \log n}{\sqrt{n \log n}}$, and suppose that $p = \ds\frac{1}{2} - \delta$ and $t = \ds\sqrt{\frac{n}{\log n}}$. Then in \textup{\textbf{Boot3($t$)}},  $A^{(4)} = A^{(3)}$ whp.
\end{lemma}

\begin{proof}
Let $n \in \N$ be sufficiently large, let $x \in Q_n$, and suppose that $x \in A^{(4)} \setminus A^{(3)}$. Then there exists a set $T \subset \Gamma(x)$ with $|T| = t$, and such that $T \subset A^{(3)} \setminus A^{(2)}$. This is because $x \in A^{(4)} \setminus A^{(3)}$ implies $|\Gamma(x) \cap A^{(3)}| \ge n/2$, and $x \notin A^{(3)}$ implies $|\Gamma(x) \cap A^{(2)}| < n/2 - t$. It is convenient to think of the vertices of $Q_n$ as subsets of $[n]$, and to assume (without loss of generality) that $x = \emptyset$, and that $T = \{\{i\} : i \in [t]\}$. We have at most $2^n$ choices for the set $T$.

Similarly, each vertex $y \in T$ must have at least $t + 1$ neighbours in $A^{(2)} \setminus A^{(1)}$. Since $x \notin A^{(2)}$, these are in $S(x,2)$, and note that each vertex of $S(x,2)$ has at most two neighbours in $T$. Thus there must exist a set $T' \subset (A^{(2)} \setminus A^{(1)}) \cap S(x,2) \cap \Gamma(T)$ with $t^2/2 \le |T'| \le t^2$. Let $|T'| = t'$. Given $T$, we have at most $t^2{{nt} \choose t^2}$ choices for $T'$.

Using the same logic one more time, each vertex $y \in T'$ must have at least $t + 1$ neighbours in $A^{(1)} \setminus A^{(0)}$, and at least $t$ of these are in $S(x,3)$. Each vertex of $S(x,3)$ has at most three neighbours in $T'$, so there must exist a set $S \subset (A^{(1)} \setminus A^{(0)}) \cap S(x,3) \cap \Gamma(T')$ with $t^3/6 \le |S| \le t^3$. Let $|S| = s$, and, considering $S$ as a 3-uniform hypergraph on $[n]$, let $\|S\| = m$. If $m < t^5$, we could use the trivial upper bound ${{nt^2} \choose s}$ on the number of choices for $S$ (given $T'$). For $m \ge t^5$ however, we shall need the following stronger bound, which follows from Lemma~\ref{count}.\\

\noindent \ul{Claim 1}: Given $T'$, $s$ and $m$, there are at most $$2^n \left( \frac{20nt}{s} \right)^{s} {{nt^2} \choose {s - 2m/n + 2t^5/n}}$$ potential sets $S$ with $|S| = s$ and $\|S\| \ge m$.

\begin{proof}[Proof of Claim 1]
Let $G$ be the 3-uniform hypergraph on $[n]$ with edge set $E(G) = \{e \in S(x,3) \cap \Gamma(T')\}$, and similarly consider $S$ to be a 3-uniform hypergraph in the obvious way, i.e., $E(S) = \{e \in S(x,3) \cap S\}$. Note that $S \subset E(G)$, $|E(G)| \le nt' \le nt^2$, $t^2 \le n \ll t^3$, and $s \le t^3$, so the result follows immediately by Lemma~\ref{count}.
\end{proof}

Now, consider the neighbourhood $D$ of $S$ in $S(x,4)$, and let $d$ be the number of edges of $Q_n$ between $S$ and $D$. Next, partition $D$ into four parts, $D_1$, $D_2$, $D_3$ and $D_4$, where each element of $D_i$ has $i$ neighbours in $S$, and let $|D_i| = d_i$. Furthermore, let $R_i = D_i \cap A^{(0)}$, and let $|R_i| = r_i$. We have
$$d \; = \; d_1 + 2d_2 + 3d_3 + 4d_4 \; = \; (n - 3)s,$$ and
$$r \; := \; r_1 + 2r_2 + 3r_3 + 4r_4 \; \ge \; s\left(\frac{n}{2} - 3t - 2\right),$$ since each vertex in $S$ has at least $n/2 - 3t $ neighbours in $A^{(0)}$, and at most two neighbours in $A^{(0)} \setminus S(x,4)$. Also,
$$\Ex(r) \; = \; pd \; = \; p(n-3)s \; \le \; s \left(\frac{n}{2} - \delta n \right),$$ and $r_i \sim \textup{Bin}(d_i,p)$ for $i = 1,2,3,4$.

In order to apply Claim 1, we shall need some bound on $\|S\|$. The following claim gives us one.\\

\noindent \ul{Claim 2}: $\|S\| \ge d_2 + 3d_3 + 6d_4$.

\begin{proof}[Proof of claim]
If $ijk\ell \in D_2$ then $ijk \in S$ and $ij\ell \in S$, say. Thus we add one pair to $d_S(i,j)$. Similarly if $ijk\ell \in D_3$, we add a pair to $d_S(i,j)$, $d_S(i,k)$ and $d_S(i,\ell)$, say, and if $ijk\ell \in D_4$, we add a pair to each of the six degrees. Finally, each pair of 3-sets are both contained in at most one 4-set.
\end{proof}

Now we apply Lemma~\ref{layer4}. Let $m = d_2 + 3d_3 + 6d_4$, so $\|S\| \ge m$. By Claim 1, we have at most
$$2^n \left( \frac{20nt}{s} \right)^{s} {{nt^2} \choose {s - 2m/n + 2t^5/n}}$$ ways of choosing the elements of $S$. Thus, recalling that $s \ge t^3/6 \gg n$, so $2^n (\delta ns)^3 \ll 2^s$, we get
\begin{align*}
& \Pr(\exists S \, | \, s, T') \; \le \; \sum_{d_1,d_2,d_3,d_4} 2^n \left( \frac{20nt}{s} \right)^{s} {{nt^2} \choose {s - 2m/n + 2t^5/n}}\\
& \hspace{7.5cm} \times \; \Pr\Big(r \ge \Ex(r) + (\delta n - 4t)s \Big)\\
& \le \; \sum_{d_1,d_2,d_3,d_4} \left( \frac{300n}{t^2} \right)^{s} \left( \frac{ent^2}{s - \frac{2m}{n} + \frac{3t^5}{n}} \right)^{s - \frac{2m}{n} + \frac{3t^5}{n}} \exp\left(- \frac{2(\delta n - 4t)^2 s^2}{D(4)} \right),
\end{align*}
where $D(4) = d_1+ 4d_2 + 9d_3 + 16d_4$, as in Lemma~\ref{layer4}. Note the replacement of $2t^5/n$ by $3t^5/n$, which is motivated by the later computation.

The rest of the proof is just a straightforward calculation. Indeed, consider
$$M(\textbf{d}) = \left( \frac{300n}{t^2} \right)^{s} \left( \frac{ent^2}{s - 2m/n + 3t^5/n} \right)^{s - \frac{2m}{n} +
\frac{3t^5}{n}} \exp\left(- \frac{2(\delta n - 4t)^2 s^2}{D(4)} \right),$$
where  $(n-3)s = d_1 + 2d_2 + 3d_3 + 4d_4$ and $m = d_2 + 3d_3 + 6d_4$, and note that
\begin{eqnarray}
D(4) & = & d_1+ 4d_2 + 9d_3 + 16d_4\nonumber\\
& = & (n-3)s + 2d_2 + 6d_3 + 12d_4 \; = \; (n-3)s + 2m.\label{triv}
\end{eqnarray}
This trivial observation will allow us to bound $M(\textbf{d})$ from above for all $(d_1, d_2, d_3, d_4)$.

Indeed, recalling that $\delta = \ds\frac{1}{2} \sqrt{\ds\frac{\log n}{n}} - \frac{\lambda \log \log n}{\sqrt{n \log n}}$ and $t = \sqrt{\ds\frac{n}{\log n}}$, first observe that
\begin{eqnarray*}
n \log n - 4\lambda n \log \log n \; \ge \; 4(\delta n - 4t)^2 & \ge & n \log n - 4\lambda n \log \log n - 32\delta tn\\
& \ge & n \log n - 3\lambda n \log \log n
\end{eqnarray*}
since $t^2 \ll n$ and $1 \le 32\delta t \ll \log \log n$. Thus, if $m \ge \alpha sn$ for some $\alpha \ge 0$, then
\begin{eqnarray}
\frac{2(\delta n - 4t)^2 s^2}{sn + 2m} & = & \frac{2(\delta n - 4t)^2 s^2}{sn} \: - \: \frac{4m(\delta n - 4t)^2 s^2}{sn(sn + 2m)}\nonumber\\[+1ex]
& \ge & \frac{s(\log n - 3\lambda \log \log n)}{2} \: - \: \frac{m(\log n - 4\lambda \log \log n)}{(1 + 2\alpha)n}.\hspace{1cm} \label{alpha}
\end{eqnarray}
Finally, we may assume that $s - 2m/n + 2t^5/n \ge 0$, as otherwise $\ds{{e(G)} \choose {s - 2m/n + 2t^5/n}} = 0$. Thus, applying (\ref{alpha}) with $\alpha = 0$, and using (\ref{triv}), we have
\begin{eqnarray*}
M(\textbf{d}) & \le & \left( \frac{300n}{t^2} \right)^{s} \left( \frac{ent^2}{t^5/n} \right)^{s - \frac{2m}{n} + \frac{3t^5}{n}} \exp\left(- \frac{2(\delta n - 4t)^2 s^2}{sn + 2m} \right)\\
& \le & \left( \frac{900n}{t^2} \right)^{s} \left( \frac{n^2}{t^3} \right)^{s - \frac{2m}{n} + \frac{3t^5}{n}} n^{-s/2 + m/n} (\log n)^{\lambda(3s/2 - 4m/n)}\\
& = &  \left( \frac{t^6}{n^3(\log n)^{4\lambda}} \right)^{m/n} \left( \frac{n^2}{t^3} \right)^{3t^5/n} \left( \frac{900n^{5/2}(\log n)^{3\lambda/2}}{t^5} \right)^s.
\end{eqnarray*}
Now, recall that $\ds\frac{n}{t^2} = \log n$, and $s \ge \ds\frac{t^3}{6}$. So if $m \le \ds\frac{sn}{4}$, then
\begin{eqnarray*}
M(\textbf{d}) & \le & \left( \frac{1}{(\log n)^{4\lambda + 3}} \right)^{m/n} n^{3t^3/\log n} \big( 900(\log n)^{3\lambda/2 + 3/2} \big)^s\\[+1ex]
& \le & \left( \frac{1}{(\log n)^{4\lambda + 3}} \right)^{s/4} e^{18s} \big( 900(\log n)^{3\lambda/2 + 3/2} \big)^s\\[+1ex]
& \le & \left( 900e^{18}(\log n)^{\lambda/2 + 3/4} \right)^s \; \le \; \big( \log n \big)^{-s/6}.
\end{eqnarray*}
However, if $m \ge \ds\frac{sn}{4}$, then we may apply (\ref{alpha}) with $\alpha = 1/4$, so
\begin{eqnarray*}
M(\textbf{d}) & \le & \left( \frac{900n}{t^2} \right)^{s} \left( \frac{n^2}{t^3} \right)^{s - \frac{2m}{n} + \frac{3t^5}{n}} \exp\left(- \frac{2(\delta n - 4t)^2 s^2}{sn + 2m} \right)\\[+1ex]
& \le & \left( \frac{900n}{t^2} \right)^{s} \left( \frac{n^2}{t^3} \right)^{s - \frac{2m}{n} + \frac{3t^5}{n}} n^{-s/2 + 2m/3n} (\log n)^{\lambda(3s/2 - 8m/3n)}\\\\[+0.5ex]
& = &  \left( \frac{t^{18}}{n^{10}(\log n)^{8\lambda}} \right)^{m/3n} \left( \frac{n^2}{t^3} \right)^{3t^5/n} \left( \frac{900n^{5/2}(\log n)^{3\lambda/2}}{t^5} \right)^s\\[+1ex]
& \le & \left( \frac{Cn^{20}(\log n)^{10\lambda}}{t^{42}} \right)^{s/12} \; = \; \left( \frac{C(\log n)^{20 + 10\lambda}}{t^2} \right)^{s/12}
\end{eqnarray*}
for some constant $C$. (Note that in both calculations we made the substitution $m = sn/4$; we could do this in the first case because $4\lambda + 3 < 0$, and in the second because $n^{10} \gg t^{18}$.) Since $t \gg \log n$, it follows again that $M(\textbf{d}) \le \big( \log n \big)^{-s/6}$ for sufficiently large $n$.

Thus, since $d_1, \ldots, d_4 \le sn$,
\begin{eqnarray*}
\Pr(\exists S \, | \, s, T') & \le & \sum_{d_1,d_2,d_3,d_4} \big( \log n \big)^{-s/6} \; \le \; \big(sn\big)^4 \big( \log n \big)^{-t^3/36},
\end{eqnarray*}
and we have at most $t^3$ choices for $s$, and at most
$2^{2n}t^2{{nt} \choose t^2}$ choices for $x$, $T$ and $T'$. Thus,
summing over all of these, we obtain
$$\Pr(\exists \, x \in A^{(4)} \setminus A^{(3)}) \; \le \; t^2 2^{2n} \left(\frac{en}{t} \right)^{t^2} \big(nt^3\big)^4 \left(\frac{1}{\log n} \right)^{t^3/36} \; = \; o(1),$$ as required.
\end{proof}

At last, we are ready to prove Theorem~\ref{sharp}.

\begin{proof}[Proof of Theorem~\ref{sharp}]
The upper bound in the theorem is exactly Corollary~\ref{upperbound}, so let $\lambda \le -2$, let $n \in \N$ be sufficiently large, let $t = \ds\frac{n}{\log n}$, let $$p(n) \; = \; \frac{1}{2} \: - \: \frac{1}{2}\sqrt{\frac{\log n}{n}} \: + \: \frac{\lambda \log \log n}{\sqrt{n \log n}},$$ and consider the \textup{\textbf{Boot3($t$)}} process on $Q_n$.

By Lemma~\ref{notall} we have $\Pr\big(A^{(3)} = V(Q_n)\big) = o(1)$, and by Lemma~\ref{boot3} we have $\Pr\big(A^{(4)} \neq A^{(3)}\big) = o(1)$. Therefore, $\Pr(A$ percolates$) = o(1)$ in the \textup{\textbf{Boot3($t$)}} process. The obvious coupling of \textbf{Boot} with \textup{\textbf{Boot3($t$)}} now shows that $\Pr(A$ percolates$) = o(1)$ in the original process also, as required.
\end{proof}

We conclude the section by briefly discussing ways in which Theorem~\ref{sharp} could be strengthened, and the limitations of our method. The alert reader will no doubt have noticed that the constant $\lambda = -2$ is not sharp; indeed, with a little more care (and no extra ideas) we could have proved that $A$ is unlikely to percolate whenever $\lambda < -3/4$.

However, our method, as it stands, cannot prove the result for any $\lambda > -1/4$. To see this, consider an ideal set $S$, with $|S| = t^3$ and $\|S\|$ small. Using our method (based on the random variable $r$), and Lemmas~\ref{normalcher} and \ref{chernoff}, it has probability at most about
$$\exp\left( - 2\delta^2 nt^3 \right)$$
of being contained in $A^{(1)} \setminus A^{(0)}$. There are about ${{nt^2} \choose {t^3}}$ such sets, and so, writing $X$ for the number of suitable sets $S \subset A^{(1)} \setminus A^{(0)}$, we get
\begin{eqnarray*}
\Pr\big( X \ge 1\big) & \le & \Ex\big(X\big) \; \approx \; {{nt^2} \choose {t^3}}\exp\left( - 2\delta^2 nt^3 \right)\\
& \approx & \left( \frac{n}{t} \exp\left( -\frac{\log n}{2} + 2\lambda \log \log n \right) \right)^{t^3} \; \approx \; \big( \log n \big)^{(2\lambda + 1/2)t^3},
\end{eqnarray*}
which is small only if $\lambda < -1/4$.

There are two obvious places in which we could potentially be leaking probability. The first is in our estimation (using $r$) of the probability that $S$ is contained in $A^{(1)} \setminus A^{(0)}$; the second is in the inequality $\Pr( X \ge 1) \le \Ex(X)$. A heuristic calculation suggests that $\textup{Var}(X)$ is not too big, and so we suspect that the first of these is in fact the problem.

Finally, we point out that two simpler changes, which one might think would improve the result, in fact do not help. Firstly, we could take the process one (or more) step(s) further (i.e., consider a \textup{\textbf{Boot4($t$)}} process), but we would just run into the same calculation, with $t^3$ replaced by $t^4$.  Alternatively, we could increase $t$; however, we cannot do so significantly, since we need the inequality $\delta t \ll \log \log n$ in order to prove inequality (\ref{alpha}).

\section{Percolation on $d$-regular graphs}\label{dregsec}

In this section we shall prove Theorem~\ref{dreg}, which uses the ideas of the previous two sections, and deduce Corollary~\ref{n^d}. Throughout this section, let $G$ be a graph as in the statement of Theorem~\ref{dreg}, so $G$ is a $d$-regular graph on $N$ vertices satisfying
$$|S(x,i) \cap \Gamma(y)| \; \le \; f_i(d)$$
for every $x,y \in V(G)$ with $y \in V(G) \setminus B(x,i-1)$, and every $i \in [k]$, where $d,N,k \in \N$, and $f_1, \ldots, f_k: \N \to \N$ are functions satisfying $$1 \: \le \: f_i(d) \: \le \: \ds\frac{\eta d}{k\log d}$$ for some small constant $\eta > 0$, to be chosen later. Since percolation on a $1$-regular graph is not very interesting, let us assume that $d \ge 2$. In fact, when $\eta$ is small, the inequality above implies that $d \ge 1/\eta$.

We begin with the upper bound; the idea is that, by Chernoff's inequality (Lemma~\ref{normalcher}), all but exponentially few of the vertices are in $A^{(1)}$, and therefore that $\Pr(x \notin A^{(k)}) \ll 1/N$. First we use Lemma~\ref{partition} to prove a version of Lemma~\ref{hyperpartition} applicable to the graph $G$.

\begin{lemma}\label{Gpartition}
Let $d,k \in \N$ and $G$ be as described above, and let $x \in V(G)$. Then there exists a partition $$S(x,k) \: = \: B_1 \cup \ldots \cup B_{m}$$ of $S(x,k)$ into $m \le d\big(f_{k-1}(d) + f_k(d)\big) + 1$ sets, such that if $y,z \in B_j$ for some $j$, then $d(y,z) \ge 3$.
\end{lemma}

\begin{proof}
Let $y \in S(x,k)$, and consider the set
$$Y \; := \; B(y,2) \cap S(x,k) \; = \; \{z \in S(x,k) : d(y,z) \le 2\}.$$ We claim that $|Y| \le d\big(f_{k-1}(d) + f_k(d)\big) + 1$. Indeed, $y$ has $d$ neighbours, of which none are in $B(x,k-2)$, at most $f_{k-1}(d)$ are in $S(x,k-1)$, and at most $f_k(d)$ are in $Y$. But if $z \notin B(x,k-1)$ then it has at most $f_k(d)$ neighbours in $S(x,k)$, by assumption. Thus
\begin{eqnarray*}
|Y| & \le & 1 \: + \: f_k(d) \: + \: df_{k-1}(d) \: + \: \big(d - f_{k-1}(d)\big)f_k(d)\\
& \le & d\big(f_{k-1}(d) + f_k(d)\big) \: + \: 1
\end{eqnarray*} as claimed.
Now, by Lemma~\ref{partition} applied to the graph $G[S(x,k)]$, it follows that the claimed partition exists.
\end{proof}

We are ready to prove the main step in the upper bound.

\begin{lemma}\label{A1toall}
Let $d,k \in \N$, $0 < \eta \le 1/6$ and $G$ be as described above. Let $x \in V(G)$ and $c > 0$. Suppose that $2\eta \le c$, and that $p$ is chosen such that, in the majority bootstrap process on $G$,
$$\Pr(y \notin A^{(1)}) < e^{-c d}.$$ Then,
$$\Pr(x \notin A^{(k+1)}) \; \le \; d^2\exp\left( -\frac{cd^k}{3^{k+1}\big(f_{k-1}(d) + f_k(d)\big)\prod_{i = 1}^{k-1} f_i(d)} \right).$$
\end{lemma}

\begin{proof}
The proof is very similar to that of Lemma~\ref{exptoall}, and so we shall give only a sketch. Let $x \in V(G)$, and let $m = d\big(f_{k-1}(d) + f_k(d)\big) + 1$. By Lemma~\ref{Gpartition}, there exist sets $B_1, \ldots, B_m$ such that $\bigcup B_j = S(x,k)$, and for each $j \in [m]$, the events $\{y \in A^{(1)}\}_{y \in B_j}$ are independent.

We proceed as in the proof of Lemma~\ref{exptoall}: suppose $x \notin A^{(k+1)}$; then for each $t \in [0,k]$ there exists a set $T(t) \subset S(x,t)$ such that $T(t) \cap A^{(k+1-t)} = \emptyset$, and
$$|T(t)| \; \ge \; \frac{d^t}{3^t \prod_{i=1}^{t-1} f_i(d)}.$$

Indeed, let $T(0) = \{x\}$, and note that since $x \notin A^{(k+1)}$, $T(0)$ satisfies the conditions. Now, suppose we have found $T(t)$ as required. Then, each $y \in T(t)$ has at most $d/2$ neighbours in $S(x,t+1) \cap A^{(k-t)}$ (since $y \notin A^{(k+1-t)}$), and at most $f_{t-1}(d) + f_t(d) \le d/6$ neighbours in $B(x,t)$ (by the properties of $G$, and since $\eta \le 1/6$), and thus it has at least $d/3$ neighbours in $S(x,t+1) \setminus A^{(k-t)}$. Moreover, each element of $S(x,t+1)$ has at most $f_t(d)$ neighbours in $S(x,t)$. Thus, by counting edges, there must exist a set $T(t+1) \subset S(x,t+1) \setminus A^{(2k-t)}$ such that $|T(t+1)| \ge d|T(t)|/3f_t(d)$, as required.

Now, define
$$N(k) \; = \; \frac{d^k}{m3^k\prod_{i = 1}^{k-1} f_i(d)},$$
and note that $|T(k)| \ge mN(k)$ and that therefore, by the pigeonhole principle, there must exist $j \in [m]$ such that $|T(k) \cap B_j| \ge \lceil N(k) \rceil \ge 1$.

Now, recall that $\Pr(y \notin A^{(1)}) < e^{-c d}$ for each $y \in B_j$, that the events $\{y \in A^{(1)}\}_{y \in B_j}$ are independent, and that $T(k) \cap A^{(1)} = \emptyset$. Observe also that, since $k\log d \le \eta d$,
$$e^{-cd} |B_j|^2 \; \le \; d^{2k}e^{-cd} \; = \; \exp\left( -cd + 2\eta d \right) \; \le \; 1.$$
Thus, by Lemma~\ref{nunlikely} we have
$$\Pr\big(|T(k) \cap B_j| \ge N(k) \big) \; \le \; 2\left( e^{-cd} \right)^{N(k)/2} \; = \; 2\exp\left( -\frac{cdN(k)}{2} \right),$$ and so
\begin{eqnarray*}
\Pr\big(x \notin A^{(2k+1)}\big) & \le & \Pr\big(\exists j\textup{ with } |T(k) \cap B_j| \ge N(k)\big) \; \le \; 2m \exp\left( -\frac{cdN(k)}{2} \right)\\
& \le & d^2\exp\left( -\frac{cd^k}{3^{k+1}\big(f_{k-1}(d) + f_k(d)\big)\prod_{i = 1}^{k-1} f_i(d)} \right)
\end{eqnarray*} as required.
\end{proof}

The upper bound in Theorem~\ref{dreg} will follow easily from Lemma~\ref{A1toall} (see `Proof of Theorem~\ref{dreg}', below), and so we now turn to the lower bound. The method is based on that of Section~\ref{lowersec}; we begin by defining the natural generalization of the \textup{\textbf{Boot1($t$)}} and \textup{\textbf{Boot3($t$)}} processes. Given $k,t \in \N$, we call the following process \textbf{Bootk($t$)}.\\[-1ex]
\begin{itemize}
\item The elements of $A^{(0)}$ are chosen independently at random, each with probability $p$.\\[-1.5ex]
\item If $0 \le m \le k-1$, then\\[-1.8ex]

$x \in A^{(m+1)}$ if $x \in A^{(m)}$ or $|\Gamma(x) \cap A^{(m)}| \ge n/2 - (k-m)t$.\\[-1.5ex]
\item If $m \ge k$, then $x \in A^{(m+1)}$ if $x \in A^{(m)}$ or $|\Gamma(x) \cap A^{(m)}| \ge n/2$.\\[-1ex]
\end{itemize}
We shall show that, if $p = \ds\frac{1}{2} - \eps$ and $t = \ds\frac{\eps d}{3k}$, then in the \textup{\textbf{Bootk($t$)}} process we have $A^{(k+1)} = A^{(k)} \neq V(G)$ with high probability. The following lemma is the key step.

\begin{lemma}\label{bootk}
Let $d,k \in \N$, $\eta > 0$ and $G$ be as described above. Let $\eps > 0$, $p = \ds\frac{1}{2} - \eps$, $t = \ds\frac{\eps d}{3k}$ and $m \in [k]$, and suppose that $12\eta \le \eps^2$. Then, in the \textup{\textbf{Bootk($t$)}} process, for every $x \in V(G)$,
$$\Pr\big(x \in A^{(m+1)} \setminus A^{(m)} \big) \; \le \; \exp\left( -\frac{\eps^{m+2} d^{m+1}}{6^{m+1}k^{m}\prod_{i=1}^{m} f_i(d)} \right).$$
\end{lemma}

\begin{proof}
Let $x \in V(G)$, and suppose that $x \in A^{(m+1)} \setminus A^{(m)}$. Then we claim that, for each $\ell \in [0,m]$ there exists a set $T(\ell) \subset S(x,\ell)$ such that $T(\ell) \subset A^{(m -\ell+1)} \setminus A^{(m-\ell)}$, and
$$\frac{t^\ell}{2^\ell\prod_{i=1}^{\ell-1} f_i(d)} \; \le \; |T(\ell)| \; \le \; t^\ell.$$
Indeed, let $T(0) = \{x\}$, let $t \in [0,m-1]$, and assume that we have found $T(\ell)$ as required. Now, let $y \in T(\ell)$, and note that, as in the proofs of Proposition~\ref{boot1} and Lemma~\ref{boot3}, since $Y \in A^{(m-\ell+1)} \setminus A^{(m-\ell)}$, it follows that $y$ has at least $t$ neighbours in $A^{(m-\ell)} \setminus A^{(m-\ell-1)}$.

Recall that $y \in S(x,\ell)$, so $y$ has at most
$$f_{\ell-1}(d) \:+\: f_\ell(d) \; \le \; \ds\frac{2\eta d}{k \log d} \; \le \; \frac{t}{2}$$ neighbours in $B(x,\ell)$, since $12\eta \le \eps\log d$. Thus $y$ has at least $t/2$ neighbours in $S(x,\ell+1) \cap A^{(m-\ell)} \setminus A^{(m-\ell-1)}$, and each vertex of $S(x,\ell+1)$ has at most $f_\ell(d)$ neighbours in $T(\ell)$, so, by counting edges,
$$|S(x,\ell+1) \cap A^{(m-\ell)} \setminus A^{(m-\ell-1)}| \; \ge \; \frac{|T(\ell)| t}{2f_\ell(d)}$$ as required.

Choose a set $T(m)$ as described, let $s = |T(m)|$, note that $T(m) \subset S(x,m) \cap A^{(1)} \setminus A^{(0)}$, and consider the neighbourhood $W$ of $T(m)$ in $S(x,m+1)$. As in the proof of Lemma~\ref{boot3}, we shall consider the number $r$ of edges between $T(m)$ and $W$ whose end-point in $W$ lies in $A^{(0)}$, and show that it is far from the expected number.

First note that each vertex in $T(m)$ has at least $\ds\frac{d}{2} - mt$ neighbours in $A^{(0)}$, and at most $f_{m-1}(d) + f_m(d)$ neighbours in $B(x,m)$, so
$$r  \; \ge \; s\left(\frac{d}{2} - mt - f_{m-1}(d) - f_m(d) \right) \; \ge \; \frac{sd}{2} \: - \: \frac{\eps sd}{2},$$ since $mt \, \le \, \ds\frac{\eps d}{3}$ and $f_i(d) \, \le \, \ds\frac{\eta d}{k\log d} \, \le \, \ds\frac{\eps d}{12}$. Next, observe that
$$\Ex(r) \; \le \; psd \; = \; \frac{sd}{2} \: - \: \eps sd.$$
Now, each vertex of $W$ has at most $f_m(d)$ neighbours in $T(m)$; let $d_i$ be the number of vertices with exactly $i$ neighbours, and recall (from Lemma~\ref{layer4}) that  $D(k) = \sum_{i=1}^k i^2 d_i$ for each $k \in \N$. By Lemma~\ref{layer4}, we have
$$\Pr\left(r \ge \Ex(r) + \frac{\eps sd}{2} \right) \; \le \; \big( \eps sd \big)^{f_m(d)} \exp\left( -\frac{(\eps sd)^2}{2D\big(f_m(d)\big)} \right).$$
Observe that $D\big(f_m(d)\big) \, \le \, f_m(d) \sum_{i=1}^{f_m(d)} i d_i \, \le \, f_m(d) sd$, note that $\ds\frac{s}{f_m(d)} \ge \ds\frac{t}{2f_m(d)} \ge \ds\frac{\eps \log d}{6\eta}$, and recall that $48\eta^2 \le \eps^3 \log d$. It then follows that $\ds\frac{\eps^2 sd}{2f_m(d)} \,\ge\, \ds\frac{\eps^3 d \log d}{12\eta} \,\ge\, 4\eta d$, and thus
\begin{eqnarray*}
\Pr\left(r \ge \Ex(r) + \frac{\eps sd}{2} \right) & \le & \exp\left(\frac{\eta d}{k \log d} \log (d^{k+1}) \right) \exp\left( -\frac{\eps^2 sd}{2f_m(d)} \right)\\[+0.5ex]
& \le & \exp\left( 2\eta d \, - \, \frac{\eps^2 sd}{2f_m(d)} \right) \; \le \; \exp\left( -\frac{\eps^2 sd}{4f_m(d)} \right)
\end{eqnarray*}
since $s \le t^m \le d^k$.

We have shown that the probability of a particular $s$-set being contained in $A^{(1)} \setminus A^{(0)}$ is small; now we simply sum over all possible $s$-sets. There are at most $\ds{{d^m} \choose s}$ choices for the set $T(m)$, and so
\begin{eqnarray*}
\Pr\big(x \in A^{(m+1)} \setminus A^{(m)} \big) & \le & {{d^m} \choose s} \exp\left( -\frac{\eps^2 sd}{4f_m(d)} \right).
\end{eqnarray*}
But $\ds\frac{\eps^2 sd}{4f_m(d)} \,\ge\, \ds\frac{\eps^2 k s\log d}{4\eta} \,\ge\, 3ms \log d$, and $\log\ds{{d^m} \choose s} \,\le\, ms \log d$, so
$${{d^m} \choose s} \exp\left( -\frac{\eps^2 sd}{4f_m(d)} \right) \; \le \; \exp\left( -\frac{\eps^2 sd}{6f_m(d)} \right).$$
Thus, using our lower bound on $s$, and recalling that $t = \ds\frac{\eps d}{3k}$, we get
\begin{eqnarray*}
\Pr\big(x \in A^{(m+1)} \setminus A^{(m)} \big) & \le & \exp\left( - \left( \frac{\eps^2 d}{6f_m(d)} \right) \left( \frac{t^m}{2^m \prod_{i=1}^{m-1} f_i(d)} \right) \right) \\
& \le & \exp\left( -\frac{\eps^{m+2} d^{m+1}}{6^{m+1}k^{m}\prod_{i=1}^m f_i(d)} \right),
\end{eqnarray*}
 as required.
\end{proof}

We are ready to prove Theorem~\ref{dreg}. The crucial property of $G= G(d)$, which we have as yet not used, will be that $N = |V(G)|$ satisfies
$$N \; \le \; \exp\left( \frac{d^k}{\big(\omega(d)k\big)^k \big(f_{k-1}(d) + f_k(d)\big) \prod_{i=1}^{k-1} f_i(d)} \right)$$ for every $d \in \N$.

\begin{proof}[Proof of Theorem~\ref{dreg}]
We shall deduce Theorem~\ref{dreg} from Lemmas~\ref{A1toall} and \ref{bootk}. In order to apply these lemmas, we need $\eta$ to be sufficiently small, so first let $\eps > 0$, recall that $f_i(d) \le f(d)$ for each $i \in [k]$, and let $d$ be sufficiently large so that
$$f(d) \; \le \; \ds\frac{\eps^2 d}{12k\log d}.$$

We begin with the upper bound. Let $p = 1/2 + \eps$, and choose the elements of $A^{(0)} \subset V(G)$ independently at random with probability $p$. Let $x \in V(G)$, and recall that $N = |V(G)|$; we shall show that
$$\Pr\left( A^{(k+1)} \neq V(G) \right) \; \le \; N\, \Pr\left(x \notin A^{(k+1)}\right) \; = \; o(1)$$ as $d \to \infty$.

Indeed, $x$ has $d$ neighbours, and so, by Lemma~\ref{normalcher},
$$\Pr\left(x \notin A^{(1)}\right) \; \le \; \Pr\big(\textup{Bin}\left(d,p \right) < d/2 \big) \; \le \; \exp\big( - \eps^2 d \big).$$
Thus, by Lemma~\ref{A1toall}, applied with $c = \eps^2$,
\begin{align*}
& N \,\Pr(x \notin A^{(k+1)}) \; \le \; N d^2 \exp\left( -\frac{\eps^2d^k}{3^{k+1}\big(f_{k-1}(d) + f_k(d)\big)\prod_{i = 1}^{k-1} f_i(d)} \right)\hspace{0.5cm}\\[+0.5ex]
& \hspace{1.9cm} \le \; \exp\left( (k\log d)^k \left( \frac{1}{(\omega(d) k)^k} - \frac{\eps^2}{3^{k+1}} \right) + 2\log d \right) \; = \; o(1)
\end{align*}
as $d \to \infty$, as required.

Now we turn to the lower bound. Let $p = 1/2 - \eps$, and again choose the elements of $A^{(0)} \subset V(G)$ independently at random with probability $p$. Let $t = \ds\frac{\eps d }{ 3k }$, and recall the \textup{\textbf{Bootk($t$)}} process, defined above. We shall show that, in this `more generous'  process, we have $A^{(k+1)} = A^{(k)} \neq V(G)$ with high probability. The result then follows by a straightforward coupling of the two processes.

The first part, that $A^{(k+1)} = A^{(k)}$ in \textup{\textbf{Bootk($t$)}}, follows immediately from Lemma~\ref{bootk}, since
\begin{eqnarray*}
\Pr\big(A^{(k+1)} \neq A^{(k)}\big) & \le & \Ex\big( |A^{(k+1)} \setminus A^{(k)}| \big) \; = \; N\,\Pr\big(x \in A^{(k+1)} \setminus A^{(k)} \big)\\
& \le & N\,\exp\left( -\frac{\eps^{k+2} d^{k+1}}{6^{k+1}k^{k}\prod_{i=1}^{k} f_i(d)} \right) \\
& \le & \exp\left( \frac{d^{k}}{k^{k}\prod_{i=1}^{k} f_i(d)} \left( \frac{1}{\omega(d)^k} - \frac{\eps^{k+2} d}{6^{k+1}} \right) \right) \; = \; o(1)
\end{eqnarray*}
as $d \to \infty$.

For the second part, that $A^{(k)} \neq V(G)$, we again use Lemma~\ref{bootk}. First note that $\Pr\left(x \in A^{(1)} \setminus A^{(0)} \right) \; \le \; \exp\big( - \eps^2 d \big)$, by Lemma~\ref{normalcher}. Now, recall that, by Lemma~\ref{bootk},
\begin{eqnarray*}
\Pr\big(x \in A^{(m+1)} \setminus A^{(m)} \big) & \le & \exp\left( -\frac{\eps^{m+2} d^{m+1}}{6^{m+1}k^{m}\prod_{i=1}^{m} f_i(d)} \right)\\[+1ex]
& \le & \exp\big( - d (\log d)^m \big) \; \le \; \exp\big( - \eps^2 d \big)
\end{eqnarray*} for every $m \in [k]$. Thus, by Markov's inequality,
\begin{eqnarray}
\Pr\left( |A^{(k)} \setminus A^{(0)}| \ge \frac{N}{4} \right) & \le & \frac{4}{N} \,\Ex\left( |A^{(k)} \setminus A^{(0)}| \right)\nonumber\\
& \le & 4 \sum_{m=0}^{k-1} \Pr\big(x \in A^{(m+1)} \setminus A^{(m)} \big)\nonumber\\
& \le & 4k \, \exp\big( - \eps^2 d \big) \; = \; o(1)\label{others}
\end{eqnarray}
as $d \to \infty$. Finally, again by Lemma~\ref{normalcher},
\begin{equation}\Pr\left(|A^{(0)}| \ge \frac{3N}{4} \right) \; \le \; \exp\left( - \frac{N}{8} \right) \; = \; o(1).\label{first}\end{equation}
Combining (\ref{others}) and (\ref{first}) gives $\Pr\big(A^{(k)} = V(G)\big) = o(1)$, as required.
\end{proof}

Before proving Corollary~\ref{n^d}, let us note that, in one direction at least, Theorem~\ref{dreg} cannot be improved substantially. Suppose $k, f(d) = O(1)$, so the bound on $N$ becomes $N = \exp\big( o(d^k) \big)$. The following example shows that the theorem is false if this bound is replaced by $N = \exp\big( O(d^{k+3}) \big)$.

\begin{eg}
Given two graphs $H_1$ and $H_2$, we write $H_1 \cup H_2$ to mean the graph with vertex set $V(H_1) \cup V(H_2)$ and edge set $E(H_1) \cup E(H_2)$, i.e., the graph obtained by putting the graphs side by side.

Let $d,k \in \N$ (with $d$ large) and $C \in \RR$, let $M = \exp\big( C d^{k+3} \big)$, let $N = d^{k+3}M$, and let $H$ be a $d$-regular graph on $d^{k+3}$ vertices which satisfies the conditions of the theorem. Then the graph
$$G = H_1 \cup \ldots \cup H_{M},$$ where the graphs $H_i$ are disjoint copies of $H$, is a $d$-regular graph on $N$ vertices which satisfies the conditions of the theorem, but
$$p_c\big(G,d/2\big) \; \ge \; 1 - \eps,$$
where $\eps = \eps(C) \to 0$ as $C \to \infty$. This follows simply because if $p < 1 - \eps$, then with high probability at least one of the sets $V(H_i) \cap A$ is empty, and so no vertex of $H_i$ is ever infected.

To show that such a graph $H$ exists, we have only to consider a random $d$-regular graph on $L = d^{k+3}$ vertices. With high probability, such a graph satisfies the conditions of Theorem~\ref{dreg}, with $f(d) = k + 3$. To see this, let $H$ be such a graph, and note that $S(x,k) \le d^k$ for each $x \in V(H)$, so (heuristically) we have
$$\Pr\big(|S(x,k) \cap \Gamma(y)| \ge k + 4 \big) \; \le \; {{d} \choose {k+4}} \left( \frac{d^k}{L} \right)^{k+4} \; \le \; \left( \frac{d^{k+1}}{L} \right)^{k+4}$$ for each $y \in V(H) \setminus B(x,k-1)$. Thus
$$\Pr\big( \exists \, x,y : |S(x,k) \cap \Gamma(y)| \ge k + 4 \big) \; \le \; L^2 \left( \frac{d^{k+1}}{L} \right)^{k+4} \; = \; o(1)$$
as required. It is straightforward to make this rough argument rigorous.
\end{eg}

Finally, let us deduce Corollary~\ref{n^d} from Theorem~\ref{dreg}. Recall that $[n]^d$ denotes the $d$-dimensional torus, i.e., the graph with vertex set $\{1, \ldots, n\}^d$ and edge set $\{xy : \ds\sum_i \big| x_i - y_i \pmod n \big| = 1\}$.

\begin{proof}[Proof of Corollary~\ref{n^d}]
Let $n = n(t)$ and $d = d(t)$ be functions satisfying the given inequalities, and let $G = [n]^d$, so $G$ is a $2d$-regular graph on $N = n^d$ vertices. Let $\omega$ be any function satisfying $1 \ll \omega(d) \ll \log \log d$. We claim that $G$ satisfies the conditions of Theorem~\ref{dreg} with $\omega(d)k = \sqrt{\ds\frac{d}{\log d}}$ and $f(d) = k + 1$.

Indeed, let $x \in V(G)$, and observe that for each $m \le k$,
$$S(x,m) \; = \; \{y \in V(G) \,:\, \sum_{i=1}^d \| x_i - y_i \|_{\Z_n} \, = \, m\},$$
where $\|x_i - y_i\|_{\Z_n}$ denotes the distance between $x_i$ and $y_i$ in $\Z_n$, the integers modulo $n$. Thus, given $y \notin B(x,m-1)$, we have
$$|S(x,m) \cap \Gamma(y)| \; = \; |\{i : x_i \neq y_i\}| \; \le \; m + 1$$ if $y \in S(x,m+1)$, and $S(x,m) \cap \Gamma(y) = \emptyset$ otherwise.

All that remains is to observe that, for some $C \in \RR$,
$$\log N \; = \; d\log n \; \le \; 2^{C\sqrt{\frac{d}{\log d}}} \; \le \; \big(\log d\big)^k \; = \; \left( \frac{d}{\big( \omega(d) k \big)^2} \right)^k$$
and hence that
$$N \; \le \; \exp\left( \frac{d^k}{\big(\omega(d)k\big)^k \big(f_{k-1}(d) + f_k(d)\big) \prod_{i=1}^{k-1} f_i(d)} \right)$$
as required.
\end{proof}

\section{Further questions and conjectures}\label{qusec}

In this section we shall briefly discuss various ways in which the work in this paper could be extended. We begin by conjecturing that the upper bound in Theorem~\ref{sharp} is sharp.

\begin{conj}
$$p_c\big(Q_n,n/2\big) \; = \; \ds\frac{1}{2} \: - \: \frac{1}{2} \sqrt{\ds\frac{\log n}{n}} \: + \: \ds\frac{ \log \log n }{ 2\sqrt{n \log n} } + o\left(\ds\frac{ \log \log n}{2\sqrt{n \log n} } \right)$$ as $n \to \infty$.
\end{conj}

Next, recall once again the result of Cerf and Manzo~\cite{CM}, that
$$p_c\left([n]^d,d\right) \: = \: o(1)$$ when $d \le \log_* n$, and Corollary~\ref{n^d}, which states that
$$p_c\left([n]^d,d\right) \: = \: \frac{1}{2} \,+\, o(1)$$ when $d \ge \eps (\log\log n)^2 \log\log\log n$. The obvious question poses itself: What happens in between?

\begin{prob}
Determine $\ds\lim_{t \to \infty} p_c\left( [n]^d, d \right)$ for every pair of functions $n = n(t)$ and $d = d(t)$ for which the limit exists.

In particular, determine the nature of the phase transition between those pairs of functions for which the limit is zero, and those for which it is non-zero.
\end{prob}

Another entirely natural question asks what happens to the critical probability if one changes the threshold function $r = r(d)$? The proof of Theorem~\ref{dreg} extends easily to the case $r = \alpha d + o(d)$ for some constant $0 < \alpha < 1$, and implies that
$$p_c\left( [n]^d, r \right) \: = \: \alpha \,+\, o(1)$$ for the same functions $n(t)$ and $d(t)$ as in Corollary~\ref{n^d}. However, the following problem is likely to be more challenging.

\begin{prob}
Determine $p_c\big([n]^d,r\big)$ for all functions $2 \le r(d) \ll d$. In particular, characterize the pairs $\big(n(d),r(d)\big)$ of functions for which $$\frac{d}{r} \, p_c\big([n]^d,r\big) \; \to \; 1$$ as $d \to \infty$.
\end{prob}

We remark that in \cite{BB}, Balogh and Bollob\'as proved that $p_c\big([2]^d,2\big)$ is very far from $\ds\frac{r}{d}$; in fact it is (up to a constant factor) equal to
$\ds\frac{2^{-2\sqrt{d}}}{d^2}$. For even sharper results on a wider class of graphs, see \cite{BBMn^d}.

Finally, there are many other $d$-regular graphs which Theorem~\ref{dreg} fails to cover; for example, those with more than $2^{2^d}$ vertices. Our final question asks, rather vaguely, for a version of the theorem which applies to such graphs.

\begin{qu}
Does there exist a set of `local' conditions which allow one to determine (or bound) the critical probability for an arbitrary $d$-regular graph on $N \ge 2^{2^d}$ vertices?
\end{qu}

Unfortunately, it appears out of reach to prove anything (in general) for fixed $d$ and $N \to \infty$.

\section{Appendix: proofs of the tools in Section~\ref{tools}}\label{appendix}

In this appendix we shall prove the simple tools used earlier. We begin by recalling Stirling's formula,
$$\sqrt{\pi n} \left( \frac{n}{e} \right)^n \; \le \; n! \; \le \; 2\sqrt{\pi n} \left( \frac{n}{e} \right)^n,$$
and by making the following basic observations.

\begin{obs}\label{facts}
\begin{enumerate}
\item[$ $]  \hspace{0cm} \\[-1ex]
\item[$(a)$] Let $x \le 1/4$. Then $e^{-x-x^2} \: \le \: 1 - x \: \le \: e^{-x}$.\\
\item[$(b)$] Let $0 \le \delta \le 1/8$. Then $e^{-2\delta - 4\delta^2} \: \le \: \ds\frac{1 - \delta}{1 + \delta} \: \le \: e^{-2\delta + 2\delta^2}$.\\
\item[$(c)$] Let $m,n \in \N$ satisfy $8m \le n$ and $4m^3 \le n^2$. Then
$$\ds{n \choose n/2 + m} \:\ge\: \ds\frac{2^{n-1}}{\sqrt{\pi n}} \, \exp\left(- \ds\frac{2m^2}{n} - 1 \right).$$
\item[$(d)$] Let $X \sim \textup{Bin}(n,p)$, where $n \in \N$ and $p \,= \,p(n) \,\ge\, \ds\frac{1}{2} \,-\, \delta$, where $\delta = \delta(n) \to 0$ as $n \to \infty$. Let $m = m(n) \in [n/2]$. Then
$$\Pr(X = m) = o\big(\Pr(X \ge m)\big)$$ as $n \to \infty$.
\end{enumerate}
\end{obs}

\begin{proof}
Part $(a)$ is straightforward. For $(b)$, note that
$$1- 2\delta \; \le \; \ds\frac{1 - \delta}{1 + \delta} \; \le \; 1 - 2\delta + 2\delta^2$$ and apply part $(a)$. For $(c)$, note that by Stirling's Formula (applied to $n!$ and $n/2!$) and part $(a)$, we have
\begin{eqnarray*}
\ds{n \choose n/2 + m} & \ge & {n \choose {n/2}} \left( \frac{n-2m}{n} \right)^m \; \ge \; \ds\frac{2^{n-1}}{\sqrt{\pi n}}  \left( 1\, -\, \frac{2m}{n} \right)^m\\[+1ex] & \ge & \ds\frac{2^{n-1}}{\sqrt{\pi n}}  \exp\left( - \frac{2m^2}{n} - \frac{4m^3}{n^2} \right)
\end{eqnarray*}
as claimed. For $(d)$, let $n \in \N$ be sufficiently large, and observe that $\Pr(X = k) = {n \choose k}p^k(1-p)^{n-k}$ for any $k \in [n]$. Thus, for any integer $0 \le t \le \sqrt{n}$, we have
\begin{eqnarray*}
\frac{\Pr(X = m + t)}{\Pr(X = m)} & \ge & \left( \frac{n-m-t}{m+t} \right)^t \left( \frac{p}{1-p} \right)^t
\; \ge \; \left( \frac{n/2-t}{n/2+t} \right)^t \left( \frac{1 - 2\delta}{1 + 2\delta} \right)^t,
\end{eqnarray*}
and so, applying part $(b)$ twice, for $\ds\frac{2t}{n}$ and for $2\delta$,
$$\frac{\Pr(X = m + t)}{\Pr(X = m)} \; \ge \; \exp\left( - \frac{4t^2}{n}  - 5\delta t - 1 \right),$$
since $t^3 \ll n^2$ and $\delta^2 t \ll \delta t$.

Finally we sum from $t = 0$ to $1/\delta$, to obtain
$$\Pr(X \ge m) \; \ge \; \Pr(X = m) \, \sum_{t=0}^{1/\delta} \exp\left( - \frac{4t^2}{n}  - 5\delta t - 1 \right)\; \ge \; \ds\frac{e^{-10}}{\delta} \; \to \; \infty$$ as $n \to \infty$.
\end{proof}

We are now ready to prove the results in Section~\ref{tools}. We begin with our reverse Chernoff bound, Lemma~\ref{chernoff}.

\begin{proof}[Proof of Lemma~\ref{chernoff}]
First let $C = 0$. We have
$$\Pr\left(S(n) \ge \frac{n}{2}\right) \; = \; \sum_{k = n/2}^n {n \choose k} p^k (1-p)^{n-k}.$$
Note that $p(1-p) = \ds\frac{1}{4} - \delta^2$, and observe that since $p \le 1/2$, the function $g(k) = {n \choose k}
p^k (1-p)^{n-k}$ is decreasing on $[n/2,n]$. Thus, for any $m \ge 0$, we have
$$\Pr\left(S(n) \ge \frac{n}{2}\right) \; \ge \; m g\left( \frac{n}{2} + m\right) \; = \; m {n \choose n/2+m} \left( \frac{1}{4} - \delta^2 \right)^{n/2} \left( \frac{p}{1-p} \right)^m.$$
Let $m = \sqrt{\ds\frac{n}{\log n}}$. Now, by
Observation~\ref{facts}$(c)$ we have
$$\ds{n \choose n/2 + m} \; \ge \; \ds\frac{2^{n-1}}{\sqrt{\pi n}} \, \exp\left( - \frac{2m^2}{n} - 1 \right),$$
by Observation~\ref{facts}$(a)$, and recalling that $8\delta^4 n \le 1$,
$$\left( \frac{1}{4} - \delta^2 \right)^{n/2} \ge 2^{-n} \exp\left( -2\delta^2 n - 8\delta^4 n \right) \; \ge \; 2^{-n} \exp\left( -2\delta^2 n - 1 \right),$$
and by Observation~\ref{facts}$(b)$,
$$\left( \frac{p}{1-p} \right)^m \; = \; \left( \frac{1 - 2\delta}{1 + 2\delta} \right)^m \; \ge \; \exp\left(-4\delta m - 16\delta^2m \right) \; \ge \; \exp\left( - 4\delta m - 1 \right).$$
Therefore,
\begin{eqnarray*}
\Pr\left(S(n) \ge \frac{n}{2}\right) & \ge & \ds\frac{m}{2\sqrt{\pi n}} \, \exp\left( - \frac{2m^2}{n} - 2\delta^2 n - 4 \delta m - 3 \right)\\
& \ge & \ds\frac{1}{2e^4 \sqrt{\pi \log n}} \, \exp\left( - 2\delta^2 n - 4 \delta \sqrt{\frac{n}{\log n}} \right),
\end{eqnarray*}
as required, since $2m^2 \le n$ and $2e^4 \sqrt{\pi} < e^{6}$. For general $C$ the proof is the same, since $C \le m/2$ for large $n$, so
$$\Pr\left(S(n) \ge \frac{n}{2} + C \right) \; \ge \; \frac{m}{2} g\left( \frac{n}{2} + m\right),$$ and $4e^4 \sqrt{\pi} < e^{6}$.
\end{proof}

Next we prove Lemma~\ref{layer4}, which generalizes Lemma~\ref{normalcher} to a weighted binomial distribution.

\begin{proof}[Proof of Lemma~\ref{layer4}]
We prove the lemma by induction on $k$. For $k=1$, it is exactly Lemma~\ref{normalcher}$(a)$, so let $k \ge 2$, and assume it is true for $k-1$. Recall that $Y_k = Y_{k-1} + kX_k$. Thus, by the induction hypothesis and Lemma~\ref{normalcher},
\begin{align*}
& \Pr\Big(Y_k \ge \Ex(Y_k) + t\Big) \; \le \; \Pr\Big(Y_{k-1} \ge \Ex(Y_{k-1}) + t \Big) \; + \; \Pr\left(X_k \ge \Ex(X_k) + \frac{t}{k} \right)\\
& \hspace{2.5cm} + \; \sum_{m=1}^{t-1} \Pr\Big(Y_{k-1} \ge \Ex(Y_{k-1}) + m \Big) \, \Pr\left(X_k \ge \Ex(X_k) + \frac{t-m}{k}\right)\\
& \hspace{0.5cm}\le \; \exp\left( - \frac{2t^2}{D(k-1)} \right) \; + \; \exp\left(- \frac{2t^2}{k^2d_k} \right) \\
& \hspace{4.5cm} + \; \sum_{m=1}^{t-1} (2t)^{k-2} \exp\left( - \frac{2m^2}{D(k-1)} - \frac{2(t-m)^2}{k^2d_k}\right).
\end{align*}
Now, simple calculus gives the maximum at $m = \ds\frac{tD(k-1)}{D(k-1) + k^2d_k}$. Thus
$$\frac{m^2}{D(k-1)} + \frac{(t-m)^2}{k^2d_k} \; \ge \; \frac{t^2 k^2 d_k D(k-1) + t^2 \big( k^2 d_k \big)^2}{k^2d_k\big(D(k-1) + k^2d_k\big)^2} \; = \; \frac{t^2}{D(k)},$$
since $D(k-1) + k^2d_k = D(k)$, and so
\begin{eqnarray*}
\Pr\Big(Y_k \ge \Ex(Y_k) + t\Big) & \le & 2\exp\left( - \frac{2t^2}{D(k)} \right) \; + \; (t-1)(2t)^{k-2} \exp\left( -\frac{2t^2}{D(k)} \right)\\
& \le & (2t)^{k-1} \exp\left( -\frac{2t^2}{D(k)} \right),
\end{eqnarray*}
and the induction step is complete.
\end{proof}

The remaining lemmas are even more straightforward.

\begin{proof}[Proof of Lemma~\ref{nunlikely}]
Recall that $pn^2 \le 1$, and that $S(n) \sim \textup{Bin}(n,p)$. Thus
$$\Pr(S(n) = m) \; \le \; {n \choose m} p^m \; \le \; (np)^m \; \le \; p^{m/2},$$ and
$$\frac{\Pr(S(n) = m+1)}{\Pr(S(n) = m)} \; \le \; \sqrt{p} \; < \; \frac{1}{2}$$ for every $m \in \N$. Therefore $\Pr(S(n) \ge m) \le 2p^{m/2}$, and the second part follows immediately.
\end{proof}

\begin{proof}[Proof of Lemma~\ref{partition}]
We apply a straightforward greedy algorithm. Taking the vertices one by one, we claim that there is some set $B_j$, such that all vertices already in $B_j$ are distance at least $k+1$ from the vertex in question. Indeed, this follows immediately from the condition
$$\left| B(x,k) \setminus \{x\} \right| \: = \: \left| \{y \in V(G) : d(x,y) \le k\} \right| \,-\, 1 \: \le \: m \,-\, 1.$$
The greedy algorithm thus gives the required partition of $V(G)$.
\end{proof}

Lemma~\ref{hyperpartition} is an immediate consequence of Lemma~\ref{partition}.

\begin{proof}[Proof of Lemma~\ref{hyperpartition}]
Let us again consider the vertices of $Q_n$ as subsets of $[n]$, and let $x = \emptyset$. Given a vertex $u \in S(x,k)$, observe that the set
$$\{v \in S(x,k) : d(u,v) \le 2k - 1\}$$ is exactly the set of $k$-subsets of $[n]$ which intersect $v$. There are at most $k{n \choose {k-1}}$ such sets, and so the result follows by Lemma~\ref{partition}.
\end{proof}

Finally, we prove the two easy lemmas.

\begin{proof}[Proof of Lemma~\ref{giveone}]
Let $S'' = S' - 1 \sim \textup{Bin}(n-1,p)$, and note that we may choose a coupling so that $S'' \le S \le S'$. Now
$$\Pr\left(S'(n) \ge m \right) \; = \; \Pr\left(S''(n) \ge m - 1 \right) \; = \;  \big( 1 + o(1) \big) \Pr\left(S''(n) \ge m \right),$$
since $\Pr\big( S''(n) = m \big) = o\big( \Pr\left(S''(n) \ge m \right) \big)$ by Observation~\ref{facts}$(d)$. Thus
$$\Pr(S(n) \ge m) \; \le \; \Pr(S'(n) \ge m) \; \le \; \big( 1 + o(1) \big) \Pr\left(S(n) \ge m \right).$$
\end{proof}

\begin{proof}[Proof of Lemma~\ref{silly}]
We have
\begin{eqnarray*}
\Pr\big(X = 1 \,|\, S(n) \ge m\big) & = & \frac{\Pr\big((X = 1) \wedge (S(n) \ge m)\big)}{\Pr\big(S(n) \ge m\big) } \\[+1ex]
& = & \frac{\Pr\big(S(n) \ge m \,|\, X = 1\big)\Pr\big(X = 1\big)}{\Pr\big(S(n) \ge m\big)},
\end{eqnarray*}
and $\Pr\big(S(n) \ge m \,|\, X = 1\big) = (1 + o(1)) \Pr\big(S(n) \ge m \big)$ by Lemma~\ref{giveone}, so the lemma follows.
\end{proof}

\section{Acknowledgements}

The authors would like to thank the Institute for Mathematical Sciences at the National University of Singapore for their hospitality during June 2006, when a large part of this research was carried out.

\end{document}